\input ./style/arxiv-general.cfg
\documentclass[aap,MSNbibl,nameyear,dvips]{arximspdf}
\makeatletter
   \@ifpackageloaded{graphicx}{}{\usepackage{graphicx}}
\makeatother
\usepackage{mathbh}

\doi{10.1214/15-AAP1107}
\volume{26}
\issue{2}
\pubyear{2016}
\firstpage{892}
\lastpage{916}
\docsubty{FLA}

\makeatletter
\newcommand{\rrvert}{\vert}
\newcommand{\llvert}{\vert}
\renewcommand{\mid}{|}
\newcommand{\binomt}[2]{{#1 \choose #2}}
\newcommand{\iint}{\int\!\!\!\int}
\newcommand{\mathds}{\mathbb}
\def\RR{\mathds{R}} 
\def\NN{\mathds{N}} 
\def\PP{\mathds{P}} 
\def\EE{\mathds{E}} 
\def\ind{\mathbh{1}} 
\def\law{\mathcal{L}} 
\def\P{\mathcal{P}} 
\def\W{\mathcal{W}} 
\def\Q{\mathcal{Q}} 
\def\N{\mathcal{N}} 
\def\i{\mathbf{i}} 
\def\F{\mathcal{F}} 
\def\G{\Pi} 
\def\eps{\varepsilon} 
\def\ELR{\mathbf{E}} 
\def\M{\mathcal{M}} 
\def\mom{M} 
\def\IP{\mathbf{P}} 
\def\C{\mathcal{C}} 
\def\L{L} 
\def\R{R} 
\def\tL{\tilde{L}} 
\def\tR{\tilde{R}} 
\def\a{\eta} 
\def\lr{\Lambda} 
\def\barlr{\bar{\Lambda}} 
\def\bi{b} 
\def\bii{b} 
\newcommand{\barcq}[1]{\bar{\alpha}_{#1}} 
\newtheorem{teo}{Theorem}
\newtheorem{lem}[teo]{Lemma}
\newtheorem{cor}[teo]{Corollary}
\newproclaim{rmk}[teo]{Remark}
\makeatother

\begin{document}
\begin{frontmatter}

\title{Quantitative propagation of chaos for generalized~Kac particle systems}
\runtitle{Generalized Kac particle systems}

\begin{aug}
\author[A]{\fnms{Roberto}~\snm{Cortez}\ead
[label=e1]{rcortez@dim.uchile.cl}\thanksref{T1}}
\and
\author[A]{\fnms{Joaquin}~\snm{Fontbona}\corref{}\ead[label=e2]{fontbona@dim.uchile.cl}\thanksref{T2}}
\runauthor{R. Cortez and J. Fontbona}
\affiliation{Universidad de Chile}
\address[A]{Department of Mathematical Engineering\\
\quad and Center for Mathematical Modeling\\
UMI(2807) UCHILE-CNRS\\
Universidad de Chile\\
Casilla 170-3, Correo 3\\
Santiago-Chile\\
\printead{e1}\\
\phantom{E-mail: }\printead*{e2}}
\end{aug}
\thankstext{T1}{Supported by Proyecto Mecesup UCH0607 Doctoral Fellowship.}
\thankstext{T2}{Supported by Fondecyt Grant 1110923, Basal-CONICYT
Center for Mathematical Modeling (CMM), and Millenium Nucleus NC120062.}

%
\received{\smonth{6} \syear{2014}}
%
\revised{\smonth{2} \syear{2015}}

%
\begin{abstract}
We study a class of one-dimensional particle systems with true (Bird
type) binary interactions,
which includes Kac's model of the Boltzmann equation and nonlinear
equations for the evolution of wealth distribution arising in kinetic
economic models. We obtain explicit rates of convergence for the
Wasserstein distance between the law of the particles and their
limiting law, which are linear in time and depend in a mild polynomial
manner on the number of particles. The proof is based on a novel
coupling between the particle system and a suitable system of
nonindependent nonlinear processes, as well as on recent sharp
estimates for empirical measures.
\end{abstract}

%
\begin{keyword}[class=AMS]
\kwd[Primary ]{60K35}
\kwd[; secondary ]{82C22}
\kwd{82C40}
\end{keyword}
\begin{keyword}
\kwd{Propagation of chaos}
\kwd{Kac equation}
\kwd{wealth distribution equations}
\kwd{stochastic particle systems}
\kwd{Wasserstein distance}
\kwd{optimal coupling}
\end{keyword}
\end{frontmatter}

\section{Introduction and main result}

\setcounter{footnote}{2}

\subsection{The kinetic equation}
We consider the collection $(P_t)_{t\geq0}$ of probability measures on
$\RR$, solution of the following nonlinear kinetic-type equation:
%
\begin{equation}
\label{eqPt} \partial_t P_t = - P_t +
\Q^+ (P_t).
\end{equation}
Here, $\Q^+$ is a generalized Wild convolution, which associates with
every measure $\mu$ on $\RR$ a new measure $\Q^+(\mu)$ given by
%
\begin{equation}
\label{eqQ+} \qquad\int\phi(u) \Q^+(\mu) (du) = \iint\frac{1}{2}\ELR\bigl(\phi(
\L u + \R v) + \phi(\tL v + \tR u)\bigr) \mu(dv) \mu(du),
\end{equation}
for all bounded measurable functions $\phi$, where $(\L,\R,\tL,\tR
)$ is
a given random vector in $\RR^4$ (with known distribution) and $\ELR$
denotes the expectation with respect to it.

Equations (\ref{eqPt})--(\ref{eqQ+}) describe the behavior of an
infinite number of objects or ``particles'' subjected to binary
interactions. The state of each particle is characterized by a scalar
$u\in\RR$, and $P_t(du)$ represents the proportion of particles in
state $u$
at time $t\geq0$. The microscopic binary interactions, which occur
randomly at constant rate, are heuristically described as follows: when
a particle at state $u$ interacts with
a particle at state $v$, their states change according to the rule
%
\begin{equation}
\label{eqinterrule} (u,v) \mapsto(\L u + \R v, \tL v + \tR u).
\end{equation}
This model is a generalization of Kac's one-dimensional simplification
of the (more realistic) Boltzmann equation for a spatially homogeneous
dilute gas in $\RR^3$, in which the interacting objects represent
actual physical particles.
Specifically, in Kac's model introduced in \citet{kac1956}, the
state of a particle is its one-dimensional velocity, and the
interactions correspond to random exchanges of velocities that occur at
binary collisions that preserve kinetic energy, so that $\L= \cos
\theta= \tL$, $\R= -\sin\theta=-\tR$, with $\theta$ randomly chosen
in $[0,2\pi)$.
We refer the reader to \citet{mischler-mouhot2013} and the
references therein for historical background on Kac and Boltzmann's equations.

A further source of models of the type described by equations (\ref
{eqPt})--(\ref{eqQ+}) is the kinetic description of the evolution of
the wealth distribution in a simplified economy, studied, for instance,
in \citet{matthes-toscani2008} (see also the references therein).
In that setting, the state of a particle represents the wealth of an
economic agent, and the binary interactions correspond to trades or
economic exchanges between them. Early versions of that model assumed
%
\begin{equation}
\label{eqLR} \llvert \L\rrvert ^p + \llvert \tR\rrvert
^p = 1\qquad\mbox{a.s.,}\qquad \llvert \tL\rrvert ^p + \llvert
\R \rrvert ^p = 1\qquad\mbox{a.s.},
\end{equation}
for some $p\geq1$ [notice that in Kac's model (\ref{eqLR}) is
satisfied with $p=2$]. In the case $p=1$, for nonnegative $\L$, $\R$,
$\tL$ and $\tR$, condition (\ref{eqLR}) can be seen as \emph{exact}
conservation of total wealth in each interaction. The weaker condition
%
\begin{equation}
\label{eqELR} \ELR\bigl(\llvert \L\rrvert ^p + \llvert \tR\rrvert
^p\bigr) = 1, \qquad\ELR\bigl(\llvert \tL\rrvert ^p +
\llvert \R \rrvert ^p\bigr) = 1,
\end{equation}
interpreted as conservation of wealth only \emph{in the mean} (so that
risky trades with possible gain or loss of total wealth in each
interaction are allowed), has also been considered in order to obtain
wider classes of equilibrium distributions for the nonlinear dynamics
[see \citet{matthes-toscani2008}, \citet
{bassetti-ladelli-matthes2011}].

%

\subsection{Particle system and propagation of chaos}
In order to rigorously justify the interpretation of the model (\ref
{eqPt})--(\ref{eqQ+}) as representing the evolution of an infinite
number of interacting particles or agents, one
considers a finite system of $N$ of such particles, which we denote
$\mathbf{X}_t = (X_t^1, \ldots, X_t^N)$, starting independently with
common law $P_0$ and such that, at each binary interaction, the states
of both involved particles are modified according to the rule (\ref
{eqinterrule}). In the terminology of particle approximations of the
Boltzmann equation, a particle system with such (true) binary
interactions is called of \emph{Bird type}, as opposed to particle
systems of \emph{Nanbu type}, in which only one particle changes its
state after interaction with some other.

Specifically, the particle system $\mathbf{X}$ has infinitesimal generator
%
\begin{equation}
\label{eqgenerator} \mathcal{A}^N \phi(\mathbf{x}) = \frac{1}{2(N-1)}
\sum_{i\neq j} \int_{\RR^4} \bigl[\phi
\bigl(\mathbf {x}+a_{ij}\bigl(\eta,x^i,x^j\bigr)
\bigr) - \phi(\mathbf{x})\bigr] \lr(d\eta)
\end{equation}
for all $\mathbf{x} = (x^1,\ldots,x^N) \in\RR^N$ and every test
function $\phi$ on $\RR^N$,
where $\a= (\xi,\zeta,\tilde{\xi},\tilde{\zeta})$ denotes\vspace*{1pt} a
generic point in $\RR^4$,
$\lr$ is the joint law of $(\L,\R,\tL,\tR)$ and
$a_{ij}(\a,u,v)$ is the vector of $\RR^N$ whose $i$th and $j$th
components are $(\xi-1)u + \zeta v$
and $(\tilde{\xi}-1)v + \tilde{\zeta}u$, respectively, and which is
equal to $0$ in the
other components.

Convergence of such a particle system, more precisely of its empirical
measures $\frac{1}{N}\sum_{i=1}^N \delta_{X^i_t}$ toward the unique
solution $(P_t)_{t\geq0}$ of the nonlinear evolution (\ref{eqPt}) as
$N$ goes to infinity, has been studied in more general frameworks and
from several points of view; see, for instance, \citet
{graham-meleard1997}, \citet{mischler-mouhot2013} and the
references therein [in particular, well posedness of (\ref{eqPt}) is
by now standard]. Since the particles are exchangeable, the convergence
of the empirical measure to $P_t$ for large $N$, as a random variable
in the space of probability measures in $\RR$ endowed with the weak
topology, is equivalent to the property of \emph{propagation of chaos}
of $\mathbf{X}_t$ with respect to $P_t$ [see \citet{sznitman1989}
for background]: for every fixed $k\in\NN$, the joint law of
$X^1_t,\ldots,X^k_t$ converges
weakly to $P_t^{\otimes k}$ as $N$ goes to $\infty$. That is, when $N$
is large, any fixed number of particles of the system behaves at time
$t$ approximately like independent random variables of law $P_t$. This
property was introduced and first established by Kac himself in
\citet{kac1956} for the particle system bearing his name, and is
nowadays known to hold for a large class of particle models, under
general mild assumptions.

\subsection{Main result}
Typically, weak convergence results are not sufficiently informative,
and one looks for more quantitative statements.
In this article, we will study the Bird-type $N$-particle system
$\mathbf{X}=(X^1,\ldots,X^N)$
and its propagation
of chaos property,
in the cases $p=1$ and $p=2$.
Our main goal is to obtain rates of convergence, as $N\to\infty$,
for the Wasserstein distance between the empirical measure of the
particle system
at time $t$ and its limiting law $P_t$,
with explicit estimates on $N$ and $t$ that grow reasonably fast as
functions of $t$.

Let $p\in\{1,2\}$ be fixed. In the case $p=2$, we will assume the
additional condition
$
\ELR(\L\R+\tL\tR) = 0$,
which is certainly satisfied in Kac's model.
As a generalization of (\ref{eqELR}), we will work under the assumption
%
\begin{equation}
\label{eqelastandinelast} \tfrac{1}{2} \ELR\bigl( \llvert \L\rrvert ^p +
\llvert \R \rrvert ^p + \llvert \tL\rrvert ^p + \llvert
\tR\rrvert ^p \bigr) \leq1.
\end{equation}
With some abuse of language, for each value of $p\in\{1,2\}$ we will
say that the model is \emph{inelastic} if the latter inequality is
strict. In that case, the interaction between particles produce an
average loss of energy when $p=2$ [see, e.g., the \mbox{inelastic} Kac
model in \citet{pulvirenti-toscani2004}] or of ``wealth'' [in the
context of \citet{matthes-toscani2008}] when $p=1$.
Also, to avoid trivial situations, in all what follows
we will assume that the model is nondegenerate, that is,
$\ELR(\llvert  \R \rrvert  +\llvert  \tR\rrvert  ) > 0$; this means
that the system produces at least some effective interactions.

Let us fix some notation.
$\P(E)$ denotes the space of probability measures on the metric space $E$.
For $\mathbf{x}\in\RR^N$ and any $i=1,\ldots,N$ we define the empirical
measures $\bar{\mathbf{x}} = \frac{1}{N} \sum_j \delta_{x^j}$
and $\bar{\mathbf{x}}^{i} = \frac{1}{N-1} \sum_{j\neq i} \delta_{x^j}$,
both being elements
of $\P(\RR)$.
Define $\mom_q(\mu) = \int\llvert  u\rrvert  ^q \mu(du)$ the absolute $q$-moment of
$\mu\in\P(\RR)$.
Given a random vector $\mathbf{Z}$ on $\RR^N$, we denote its law by
$\law(\mathbf{Z}) \in\P(\RR^N)$,
and the joint law of its first $k$ components by $\law^k(\mathbf{Z})
\in\P(\RR^k)$.

Recall that for $\mu,\nu\in\P(\RR^k)$
their $p$-Wasserstein distance $\W_p(\mu,\nu)$ is defined to be the
cost of the optimal transfer plan between $\mu$ and $\nu$, that is,
\[
\W_p(\mu,\nu) = \biggl(\inf_{\pi} \int
_{\RR^k \times\RR^k} d_{k,p}(\mathbf {x},\mathbf{y})^p
\pi(d\mathbf{x},d\mathbf{y}) \biggr)^{1/p} = \Bigl(\inf
_{\theta,\vartheta} \EE d_{k,p}(\theta,\vartheta)^p
\Bigr)^{1/p},
\]
where the first infimum is taken over all measures $\pi$ on $\RR^k
\times\RR^k$
with marginals $\mu$ and $\nu$, and the second infimum is taken over
all pairs of
random vectors $\theta$ and $\vartheta$ such that $\law(\theta) =
\mu$
and $\law(\vartheta) = \nu$ [see, e.g., \citet
{villani2009} for background on Wasserstein distances].
We will use the \emph{normalized} distance $d_{k,p}$ on $\RR^k$
given by
%
\begin{equation}
\label{eqdkp} d_{k,p}(\mathbf{x},\mathbf{y}) = \Biggl(
\frac{1}{k} \sum_{i=1}^k \bigl
\llvert x^i - y^i\bigr\rrvert ^p
\Biggr)^{1/p},
\end{equation}
which is natural when one cares about the dependence on the dimension.

In order to obtain good rates of convergence in $N$ which moreover are
well behaved with respect to $t$, it is
convenient to introduce the concave function
\[
\alpha_{q} = 1 - \tfrac{1}{2} \ELR\bigl(\llvert \L\rrvert
^q + \llvert \R \rrvert ^q + \llvert \tL\rrvert
^q + \llvert \tR\rrvert ^q\bigr)\qquad\forall q\geq0.
\]
We also define
\[
q^* = \sup\bigl\{q\dvtx  \mom_q(P_0)<\infty,
\alpha_{q}>0\bigr\}.
\]
%
These objects play an important role in \citet{matthes-toscani2008},
since when $p=1$ and $q^*$ is nontrivial (i.e., $1<q^*<\infty$),
$q^*$ corresponds to the Pareto index of the stationary distribution of $P_t$.
More importantly, in the present context,
the moments of order $q<q^*$ of $P_t$ can be controlled uniformly in
time (see Lemma~\ref{le5} below).
Assuming (\ref{eqelastandinelast}) and $\mom_p(P_0)<\infty$, the
concavity of $\alpha_{q}$ implies that either $q^* \in[p,\infty]$
or $q^* = -\infty$. Also, define for all $q \in\{p\} \cup(p,q^*)$
\[
\barcq{p,q} = \inf_{p\leq r \leq q} \alpha_{r} = \min(\alpha
_{p},\alpha_{q}).
\]
Note that if $\alpha_{p} = 0$, then $\barcq{p,q} = 0$ for all such $q$,
so this function is meaningful only in the case $\alpha_{p} > 0$, in which
case it will be useful to obtain uniform (in time) estimates.

We are now ready to state our main theorem (see also Corollary \ref
{cortrajectorial}
for a trajectorial result).

%
\begin{teo}
\label{teomain2}
Let $(P_t)_{t\geq0}$ be the unique solution of (\ref{eqPt})
and let $\mathbf{X}$ be the particle system starting with law
$P_0^{\otimes N}$ and with generator (\ref{eqgenerator}).
For $p=1$ or $p=2$,
assume $\alpha_{p}\geq0$ and $\mom_p(P_0) < \infty$.
If $p=2$, assume also that $\ELR(\L\R+ \tL\tR) = 0$ and $q^*>2$.
Then:
\begin{itemize}
\item for any $q\in\{1\}\cup(1,q^*)$ and
any $\gamma< (2+1/q)^{-1}$ in the case $p=1$, or
\item for any $q\in(2,q^*)$, $q\neq4$ and for $\gamma= \min
(1/3,\frac
{q-2}{2q-2})$
in the case $p=2$,
\end{itemize}
there exists a constant $C$,
depending on $p$, $q$, $\gamma$ and some moments of $P_0$ and $(\L,\R,\tL,\tR)$
of order at most $q$, such that:
\begin{longlist}[(ii)]
\item[(i)] for all $k \leq N$ and for all $t\geq0$,
\[
\W_p^p\bigl(\law^k(\mathbf{X}_t),P_t^{\otimes k}
\bigr) \leq C \biggl( \frac{t(1+t)^{p-1}e^{-(p/q)\barcq{p,q}
t}}{N^\gamma
} + \ind_{k\neq1}\frac{k \min(1,t) e^{-\alpha_{p} t}}{N}
\biggr),
\]
\item[(ii)] for all $t\geq0$,
\[
\EE\W_p^p(\bar{\mathbf{X}}_t,P_t)
\leq\frac{C (1+t)^p e^{-(p/q)\barcq{p,q} t} }{N^\gamma}.
\]
\end{longlist}
\end{teo}

%
\begin{rmk}
$\bullet $ The power $\gamma$ in Theorem \ref{teomain2} is a consequence of
using recently established sharp quantitative estimates in Wasserstein
distance for the empirical measures of exchangeable or i.i.d. collections of random variables [which improve or extend a classical
result in \citet{rachev-ruschendorf1998}].
More specifically,
the rate $N^{-\gamma}$ with $\gamma< (2+1/q)^{-1}$ in the case $p=1$
comes from
Theorem 1.2 of \citet{hauray-mischler2014}, whereas the value
$\gamma= \min(1/3,\frac{q-2}{2q-2})$
in the case $p=2$ comes from Theorem 1 of \citet{fournier-guillin2013}.
On the other hand, the dependence on $t$ results from our estimates,
which rely on Gronwall's lemma.

$\bullet$  The restriction $q\neq4$ in the case $p=2$ comes from Theorem~1
of \citet{fournier-guillin2013}.
As those authors mention, the case $q=4$ would produce additional
logarithmic terms, which in
our case translate into a rate of order $N^{-1/3}$ times a logarithmic
function of $N$.

$\bullet$ In the elastic case (i.e., $\alpha_{p} = 0 = \barcq{p,q}$),
(i)  and (ii)
give estimates that grow linearly with time (in the case $p=2$ both sides
are squared).
In the inelastic case, which corresponds to $\alpha_{p}, \barcq{p,q} >0$,
all estimates are uniform in time.


$\bullet$ From a physical point of view, it is interesting to consider
models where infinitely many particles interact over finite time intervals,
such as the Kac equation without cutoff.
The techniques used in the proof of Theorem \ref{teomain2} can also be
applied to cutoffed approximations of that equation and, in the case
that [in the notation of \citet{desvillettes-graham-meleard1999}]
the classical condition $\int_0^\pi\theta^2 \beta(\theta)\,d\theta
<\infty$ on the cross-section function $\beta\dvtx [-\pi,\pi]\to\RR_+$ is
satisfied, they yield a constant that does not
depend on the cutoff parameter; see Remark \ref{rmknon-cutoff}.
\end{rmk}

\subsection{Particular cases and comparison with known results}

\subsubsection{The Kac equation}
Note that if the stronger condition (\ref{eqLR}) is satisfied (or
holds with $\leq$ instead of equality),
then $\llvert  \L\rrvert  $, $\llvert  \R \rrvert  $, $\llvert  \tL\rrvert  $ and $\llvert  \tR\rrvert  $ are all $\leq1$ a.s., which
implies that $\alpha_{q}$ is strictly increasing with $q$. Thus,
$\barcq
{p,q} = \alpha_{p}$ for $q\geq p$
and the value of $q^*$ will depend only on the finiteness of the
moments of $P_0$.
In Kac's model, since (\ref{eqLR}) is satisfied for $p=2$, if $P_0$
has finite moment of order $4+\eps$, then $q^*>4$ and Theorem \ref
{teomain2} gives
\[
\EE\W_2^2(\bar{\mathbf{X}}_t,P_t)
\leq\frac{C (1+t)^2}{ N^{1/3}}.
\]
Several similar results can be found in the literature.
The closest one corresponds to quantitative rates for the Nanbu system
associated with Kac's model, which are found, for instance, in the proof
of Proposition 6.2 of \citet{fournier-godinho2014}. The authors
state there a $\W_2^2$ convergence rate that also depends quadratically
on $t$ and is optimal on $N$, in the sense that it is equal to the $\W
_2^2$ rate of convergence of the empirical measure of an i.i.d. sample
toward their common law. The latter is of order $N^{-1/2}$, according
to Theorem~1 of \citet{fournier-guillin2013}. Thus, the Bird-type
particle system seems to produce a slower rate of convergence than the
corresponding Nanbu-type system.
An interesting question is whether
this difference is a mere consequence of the techniques used in our proof
(more specifically, some order is lost when one uses Lemma \ref{lemW2Ymu})
or is intrinsically related to the type of binary
interactions (Bird or Nanbu) in the system.

A similar result as the one of \citet{fournier-godinho2014} can be
found in \citet{fournier-mischler2013} where, motivated by the
numerical approximation of the Boltzmann equation for hard spheres,
hard potentials and Maxwellian gases, a~pathwise coupling argument was
developed for Nanbu particle systems, which extends a coupling
construction based on optimal transport developed in \citet
{fontbona-guerin-meleard2009}. That pathwise approach, however, does
not readily extend to the particle systems of Bird type we are
interested in, which in turn provide a physically more transparent
description of the relevant interaction phenomena.

As for the Bird particle system, in \citet{graham-meleard1997} the
authors obtain an explicit rate in
total variation distance on the path space, between the law of one particle
and the law of the nonlinear process (to be introduced later).
However, due to the generality of their hypotheses and the strong
pathwise distance they use, the convergence rate depends exponentially
on the length of the time interval that is considered.
Similarly, in Theorem 4.3 of \citet{desvillettes-graham-meleard1999}
the authors state a propagation of chaos result in $\W_2$ for the law
at time $t$ of one particle in the system with
cutoff, toward the law $P_t$ of the nonlinear dynamics without cutoff.
Since some relations between
$N$ and the cutoff parameter must be satisfied when removing the
latter, that result gives estimates that are logarithmic in $N$ and grow
exponentially with $t$.

On the other hand,
the general theory developed in \citet{mischler-mouhot2013}
provides a framework and a methodology to establish quantitative (in
$t$ and $N$) propagation of chaos estimates which can be applied in the
present framework. For instance, in their Theorem 5.2, a $\W_1$
estimate for the Boltzmann equation in the Maxwell molecules case is
obtained, which is uniform in time and decays with $N$ in a polynomial
way [see also step~3 of the proof of Theorem 8 in \citet
{carrapatoso2014B} for results in $\W_2$ distance]; we expect that
similar bounds can be obtained with their techniques for the Kac model.
The actual dependence on $N$ they give seems however hard to trace in
general, and we have not been able to deduce with their techniques an
estimate in Wasserstein distance as sharp as ours in terms of~$N$.
Also, their
approach does not provide any information on the way in which
trajectories of particles get closer to those of the limiting
processes. On the other side, unfortunately our techniques (ultimately
relying on Gronwall's lemma) do not seem to yield uniform in time
estimates for the elastic Kac equation,
even if $P_0$ were compactly supported.

Finally, we observe that for the inelastic Kac model, $\barcq{2,q} =
\alpha_{2} >0$ for all $q$, hence Theorem \ref{teomain2} does give a
uniform-in-time rate in that case:
\[
\EE\W_2^2(\bar{\mathbf{X}}_t,P_t)
\leq\frac{C e^{-(\alpha _{2}/(2+\eps))t}}{ N^{1/3}},
\]
if $P_0$ is again assumed to have finite moment of order $q>4$.
This exponential decay is not surprising: when $P_0$ has finite moment of
order 2 it is known that $M_2(P_t)$ decays exponentially fast; see
\citet{pulvirenti-toscani2004}.
Nevertheless, to our knowledge, our quantitative (in $N$) propagation
of chaos result is new for the inelastic Kac model.

\subsubsection{Models for economic exchanges and wealth distribution}
Working with $p=1$ as in \citet{matthes-toscani2008} and assuming
only that the first moment of $P_0$ is finite, Theorem \ref{teomain2}
gives in the elastic case
a propagation
of chaos result for $\W_1$ of order almost $N^{-1/3}$, with estimates
growing linearly with time.
Any additional finite $q$-moment
of $P_0$ (with $q<q^*$) can be used to improve the rate in $N$, up to
almost $N^{-1/(2+1/q^*)}$.
In the case of exact conservation of wealth [condition~(\ref{eqLR})]
we have $q^*=\infty$
and we obtain a rate of $N^{-(1/2 -\eps)}$, which is almost optimal
according to Theorem 1 of \citet{fournier-guillin2013}.
To our knowledge, this is the first quantitative propagation of chaos result
for kinetic equations modeling the evolution of wealth distribution.

\subsection{The nonlinear process and idea of the proof}

Following ideas pioneered by Tanaka in the case of the Boltzmann
equation [see \citet{tanaka1978} and \citet{tanaka1979}], it
is also possible to establish the convergence of the \emph{pathwise}
law of a particle, to the law
of some process obtained by the following construction:
consider a Poisson point measure $\M$ on
$\RR_+ \times\RR^2 \times\RR$ with\vspace*{2pt} intensity $dt \barlr(d\xi,d\zeta) P_t(dv)$,
where $\barlr= \frac{1}{2}(\law(\L,\R) + \law(\tL,\tR))$,
and let $(V_t)_{t\geq0}$ be the jump process on $\RR$ defined as the unique
solution starting with law $P_0$ of the stochastic equation
%
\begin{equation}
\label{eqnlpr} dV_t = \int_{\RR^2} \int
_{\RR} \bigl[(\xi-1) V_{t^-} + \zeta v\bigr] \M
(dt,d\xi,d\zeta,dv).
\end{equation}
It is not hard to see that such a jump process $V$ exists, it is
uniquely defined
and it satisfies $\law(V_t) = P_t$ for all $t\geq0$.
We call $P$ the pathwise law of $V$, and any process with law $P$ is
called a \emph{nonlinear process};
it represents the trajectory of any fixed particle
in the (infinite) population subjected to the random interactions
described above in (\ref{eqinterrule}).

To\vspace*{1pt} prove our results, we will couple the Bird particle system
$\mathbf{X}_t$ with a system $\mathbf{U}_t = (U_t^1,\ldots,U_t^N)$
where each $U^i$ is
a copy of the nonlinear process $V$, constructed in such a way that it
remains close to $X^i$.
To achieve this, we will use techniques of optimal coupling inspired by
those used
in \citet{fontbona-guerin-meleard2009} and \citet
{fournier-mischler2013}, in order to carefully choose the jumps of the
nonlinear process $U^i$ as similar as possible to those of the
particle $X^i$. However, contrary to those papers which deal with
Nanbu-type particle systems (in which each randomness source acts on
the trajectory of only one of the particles), ensuring closeness of
$X^i$ and $U^i$
simultaneously for all $i=1,\ldots,N$ will imply
that the processes $U^1,\ldots,U^N$ are \emph{not} independent.
Therefore, to obtain the desired estimates we will need, in a second
step, to ``decouple''
the system $\mathbf{U}_t$ as $N$ goes to infinity, which we will be
able to do with estimates that
are uniform in time; see Lemma \ref{lemdecoupling} below.

Let us point out that the coupling construction we will introduce can
in principle be replicated in
higher dimensions, and with more general interaction rules, which is
why we preferred to avoid the use of specific one-dimensional features
in its construction; see, for instance, Remark \ref{rmkpostthmcoupling}. We thus expect these techniques
to be applicable in physically more relevant situations, hopefully
including (at least some instances of) the Boltzmann equation. Also, we
think it should be possible to adapt this coupling construction in
order to quantitatively study ``Bird-type'' Brownian particle
approximations of a certain Gaussian white-noise driven nonlinear
process, associated with the Landau equation arising in the grazing
collisions limit of the Boltzmann equation. Such a process was studied
in \citet{funaki1984} and \citet{guerin2004}, and a particle
approximation result with a ``Nanbu type'' Brownian particle system was
proved in \citet{fontbona-guerin-meleard2009}, by means of a
coupling construction based on optimal transport. The corresponding
particle system of Bird-type is studied in \citet{carrapatoso2014}
using the functional tools developed in \citet
{mischler-mouhot2013}, but there seems to be so far no suitable
coupling argument available in order to deal with such class of
particle systems.\footnote{When the present work was just finished, the
authors learned from Nicolas Fournier that the latter question was
currently being studied by him, Fran\c cois Bolley and Arnaud Guillin.}

\subsection{Plan of the paper}
In Section~\ref{secprel}, we give the explicit construction of the
particle system $\mathbf{X}_t$, and more importantly, we couple it with
the system
$\mathbf{U}_t = (U_t^1,\ldots,U_t^N)$ of dependent nonlinear processes
that we will use
throughout the rest of this article. In Section~\ref{secscheme}, we
prove Theorem \ref{teomain2}.
The proof of some intermediate lemmas, including statements of
Section~\ref{secprel} and the ``decoupling'' of the process $\mathbf{U}_t$,
is left for the final Section~\ref{secproofs}.

\section{Coupling of the particle system and the nonlinear processes}
\label{secprel}

\subsection{The particle system}
Let us fix the number of particles $N\in\NN$.
Although most of the subsequent objects will depend on $N$, for
notational simplicity we will not
make this dependence explicit.
We will define both the particle system $\mathbf{X}$ and the nonlinear
processes $\mathbf{U}$
by means of integral equations driven by the same
Poisson point measure. To this end, let us first
introduce the function $\i\dvtx [0,N) \to\{1,\ldots,N\}$ given by
$\i(\rho) = \lfloor\rho\rfloor+ 1$, and the set $\C\subseteq[0,N)^2$
\[
\C= \bigl\{ (\rho,\sigma) \in[0,N)^2\dvtx  \i(\rho) \neq\i(\sigma) \bigr
\}.
\]
Note that $\llvert  \C\rrvert   = N(N-1)$.
As in (\ref{eqgenerator}), denote $\a= (\xi,\zeta,\tilde{\xi
},\tilde{\zeta})$ a generic
point in $\RR^4$
and $\lr= \law(\L,\R,\tL,\tR)$.
Now, let $\N(dt,d\a,d\rho,d\sigma)$ be a Poisson point measure on
$[0,\infty)\times\RR^4 \times[0,N)^2$
with intensity
\[
\frac{N}{2} \,dt \lr(d\a) \,d\rho \,d\sigma\frac{1}{\llvert  \C\rrvert  }
\ind_{\C
}(\rho,\sigma) = \frac{1}{2(N-1)} \,dt \lr(d\a) \,d\rho \,d\sigma
\ind_{\C}(\rho,\sigma).
\]
In words, $\N$ picks atoms in $[0,\infty)$ at constant rate of $N/2$, and
for each such atom it also independently samples a tuple $(\xi,\zeta,\tilde{\xi},\tilde{\zeta}
)$ from $\lr$
and a pair $(\rho,\sigma)$ uniformly on $\C$.
We will use $(\rho,\sigma)$
to choose the indices of the particles that interact at each jump.
Consider also $N$ independent random variables $(X_0^1, \ldots, X_0^N)
=: \mathbf{X}_0$, independent
from $\N$, each having distribution $P_0$.
Finally, set $\F= (\F_t)_{t\geq0}$ to be the complete right
continuous filtration generated by $\mathbf{X}_0$ and $\N$. We denote
$\PP$ and $\EE$ the probability and expectation in the corresponding
probability space.

The particle system $\mathbf{X} = (X^1,\ldots,X^N)$ is defined as the
solution, starting from~$\mathbf{X}_0$, of the following integral equation:
%
\begin{equation}
\label{eqXt} \qquad d\mathbf{X}_t = \int_{\RR^4} \int
_{[0,N)^2} \sum_{i,j=1}^N
\ind _{\{\i
(\rho) = i, \i(\sigma) = j\}} a_{ij} \bigl(\a, X_{t^-}^i,
X_{t^-}^j\bigr) \N(dt,d\a,d\rho,d\sigma).
\end{equation}
[Recall that $a_{ij}(\a,u,v)$ is the vector of $\RR^N$ whose $i$th and
$j$th components are $(\xi-1)u + \zeta v$
and $(\tilde{\xi}-1)v + \tilde{\zeta}u$, resp., and is equal
to $0$ in the other
components].
Given the timely ordered atoms $(t_n,\a_n,\rho_n,\sigma_n)_{n\geq0}$
of $\N$ (i.e., $t_n\leq t_{n+1}$ for all $n\geq0$),
a solution of this equation can be constructed as follows:
recursively define $\mathbf{X}_{t_n}$ as
%
\begin{equation}
\label{eqXtn} X_{t_n}^\ell= \cases{ \xi_n
X_{t_{n-1}}^{i} + \zeta_n X_{t_{n-1}}^{j},
&\quad$\ell =i$,
\vspace*{3pt}\cr
\tilde{\xi}_n X_{t_{n-1}}^{j} +
\tilde{\zeta}_n X_{t_{n-1}}^{i}, &\quad$\ell=j$,
\vspace*{3pt}\cr
X_{t_{n-1}}^\ell, &\quad$\ell\neq i,j$,}
\end{equation}
where $(i,j) = (\i(\rho_n),\i(\sigma_n))$, and set $\mathbf{X}_t =
\mathbf{X}_{t_n}$ for all $t\in(t_n,t_{n+1})$.
Uniqueness for (\ref{eqXt}) also holds, since there is
no choice to make in this construction. It is straightforward to verify that
$\mathbf{X}$ has generator (\ref{eqgenerator}).

Thus, the system $\mathbf{X}$ is what we want it to be: at rate $N/2$
we choose two distinct indices $i=\i(\rho)$
and $j=\i(\sigma)$, and then we update the particles $X^i$ and $X^j$
according to the rule
described in (\ref{eqinterrule}). The fact that we use continuous
variables $(\rho,\sigma)$ to
choose the indices $(i,j)$ (instead of a discrete pair chosen uniformly from
the set $\{1,\ldots,N\}^2\setminus\{i=j\}$) will be crucial to define our
system $\mathbf{U}$ of $N$ nonlinear processes.

\subsection{Coupling with the nonlinear processes}
From (\ref{eqXt}), it follows that for any $i = 1,\ldots,N$, the
process $X^i$ satisfies
%
\begin{equation}
\label{eqXtk} dX_t^i = \int_{\RR^2}
\int_{[0,N)} \bigl[(\xi-1)X_{t^-}^i +
\zeta X_{t^-}^{\i
(\tau)}\bigr] \N ^i(dt,d\xi,d\zeta,d
\tau),
\end{equation}
where $\N^i$ is defined as
%
\begin{eqnarray}
\label{eqNk} \N^i(dt,d\xi,d\zeta,d\tau) &=& \N\bigl(dt,\bigl(d\xi
\times d\zeta\times\RR^2\bigr),[i-1,i),d\tau\bigr)
\nonumber\\[-8pt]\\[-8pt]\nonumber
&&{} + \N\bigl(dt,\bigl(\RR^2 \times d\xi\times d\zeta\bigr),d
\tau,[i-1,i)\bigr).
\end{eqnarray}
Clearly, $\N^i$ is a Poisson point measure on $[0,\infty)\times\RR
^2\times[0,N)$ with intensity
\[
dt \barlr(d\xi, d\zeta) \frac{ d\tau}{N-1} \ind_{A^i}(\tau),
\]
where $\barlr= \frac{1}{2}(\law(\L,\R) + \law(\tL,\tR))$,
and $A^i = [0,N)\setminus[i-1,i)$. In other words, $\N^i$ selects only
the atoms
of $\N$ that produce a jump of $X^i$, that is, the atoms in which $\i
(\rho) = i$ or $\i(\sigma) =i$.

Let\vspace*{1pt} us examine the expression (\ref{eqXtk}) in more detail. First,
note that since $\tau$ is chosen uniformly
in $A^i$, the variable $X_{t^-}^{\i(\tau)}$ corresponds to a sample
from the (random)
probability measure
$\bar{\mathbf{X}}^{i}_{t^-} = \frac{1}{N-1} \sum_{j\neq i} \delta
_{X_{t^-}^j}$.
Thus,\vspace*{-2pt} from the point of view of the process $X^i$, the dynamics
is as follows: at rate $1$, a number $v = X_{t^-}^{\i(\tau)}$ is
sampled from the measure $\bar{\mathbf{X}}^{i}_{t^-}$,
and then the value of the process is updated according to the rule
$X_{t^-}^i \mapsto\xi X_{t^-}^i + \zeta v$,
where $(\xi,\zeta)$ is chosen with law $\barlr$.

Comparing (\ref{eqnlpr}) and (\ref{eqXtk}),
the key observation is the following: if for each jump time $t$ one
replaces $X_{t^-}^{\i(\tau)}$ in (\ref{eqXtk})
with a realization $v$ of the law $P_t(dv)$, the resulting process has
law $P$.
In view of this, we would like to define the system of nonlinear
processes $\mathbf{U}=(U^1,\ldots,U^N)$ based on this idea,
but using a realization of $P_t$ that is optimally coupled to
the realization
$X_{t^-}^{\i(\tau)}$ of the measure $\bar{\mathbf{X}}^{i}_{t^-}$.
In doing this, some measurability issues need to be taken into account.

%
\begin{lem}[(Coupling)] \label{lemcoupling}
For every $p\geq1$ and $i\in\{1,\ldots,N\}$
there exist a measurable mapping $\G^i\dvtx \RR_+ \times\RR^N \times A^i
\to\RR$,
$(t,\mathbf{x},\tau)\mapsto\G_t^i(\mathbf{x},\tau)$,
with the following property: for every $t\geq0$ and $\mathbf{x}\in
\RR^N$,
if $\tau$ is\vspace*{1pt} uniformly chosen from $A^i$, then the pair $(\G
_t^i(\mathbf{x},\tau),x^{\i(\tau)})$
is an optimal coupling between $P_t$ and $\bar{\mathbf{x}}^{i} =
\frac
{1}{N-1} \sum_{j\neq i} \delta_{x^j}$
with respect to the cost function $c(u,v) = \llvert  u-v\rrvert  ^p$.
Moreover, if $\mathbf{Y}$ is any exchangeable random vector in $\RR
^N$, then
$\EE\int_{j-1}^j \phi(\G_t^i(\mathbf{Y},\tau)) \,d\tau= \langle P_t,
\phi\rangle$
for any $j\in\{1,\ldots,N\}$,
$j\neq i$, and any bounded measurable function $\phi$.
\end{lem}

For simplicity, in our notation we have not made explicit the
dependence of $\G_t^i$ on $p$ (however, see Remark \ref
{rmkpostthmcoupling}).
Now, we can define $U^i$ as the solution of
%
\begin{equation}
\label{eqUtk} dU_t^i = \int_{\RR^2}
\int_{[0,N)} \bigl[(\xi-1)U_{t^-}^i +
\zeta\G _t^i(\mathbf {X}_{t^-},\tau)\bigr]
\N^i(dt,d\xi,d\zeta,d\tau),
\end{equation}
where $\N^i$ is the same Poisson point measure as in (\ref{eqXtk}).
The proof of Lemma \ref{lemcoupling} will imply that the mapping
$((t,\omega),\xi,\zeta,\tau)\mapsto(\xi-1)U_{t^-}^i (\omega)+
\zeta\G
_t^i(\mathbf{X}_{t^-}(\omega),\tau)$ above is measurable with respect
to the product of the predictable sigma field [in~$(t,\omega)$] and the
Borel sigma field of $\RR^2\times[0,N)$. This ensures that the
integral in (\ref{eqUtk}) has the usual properties of integrals with
respect to Poisson point processes.

We summarize our construction in the following.

%
\begin{lem} \label{lemXU}
Let $p\geq1$ be fixed. For each $i=1,\ldots,N$ there is a unique
solution $U^i$ of (\ref{eqUtk}),
and it is a nonlinear process.
Moreover, the collection $(X^1,U^1),\ldots,(X^N,U^N)$ is exchangeable.
\end{lem}

Thus, the system $\mathbf{U} = (U^1,\ldots,U^N)$
is indeed a tuple of $N$ nonlinear processes. However, as we already
mentioned, they are \emph{not
independent}, since $\N^i$ and $\N^j$ share a portion of $\N$,
namely, the atoms of $\N$ whose coordinates $(\rho,\sigma)$ lie in
$[i-1,i)\times[j-1,j)$ or $[j-1,j)\times[i-1,i)$.
In particular, whenever such an atom occurs\vspace*{1pt} the processes
$U^i$ and $U^j$ jump simultaneously, using a single realization of $(\L,\R,\tL,\tR)$,
and samples of $P_t$ that also are correlated.

\section{Proof of the main result}
\label{secscheme}

Before proving our results, let us first state two lemmas
that constitute our basic tools; they will be proven in Section~\ref
{secproofs}.
The first one provides uniform bounds for
the moments of $P_t$; it can be seen as a version
of Theorem 3.2 in \citet{matthes-toscani2008}.

%
\begin{lem}[(Moment bounds)]
\label{lemmomentbounds}\label{le5}
For $p=1$ or $p=2$, assume $\alpha_{p} \geq0$ and $\mom_p(P_0) <
\infty$.
If $p=2$, assume also that $\ELR(\L\R+ \tL\tR) = 0$.
Then for any $q\in\{p\}\cup(p,q^*)$ there exists
a constant $C$, depending on $q$ and some moments of $P_0$ and $(\L,\R,\tL,\tR)$
of order at most $q$, such that
\[
\mom_q(P_t) \leq C e^{- \barcq{p,q} t }\qquad\forall t
\geq0.
\]
\end{lem}

The second lemma is fundamental in our developments since it decouples
the nonindependent nonlinear processes uniformly in time, even in
the case $\alpha_{p} = 0$:

%
\begin{lem}[(Decoupling)] \label{lemdecoupling}
For $p=1$ or $p=2$,
assume $\alpha_{p}\geq0$ and $\mom_p(P_0) < \infty$.
If $p=2$, assume also that $\ELR(\L\R+ \tL\tR) = 0$.
Then there exists a constant~$C$, depending only on the $p$-moment
of $P_0$ and $(\L,\R,\tL,\tR)$, such that for all $k=2,\ldots,N$ and
$t\geq0$,
\[
\W_p^p\bigl( \law^k(\mathbf{U}_t),
P_t^{\otimes k}\bigr) \leq\frac{C (k-1) \min(1,t) e^{-\alpha_{p} t}}{N-1}.
\]
\end{lem}

\begin{pf*}{Proof of Theorem \ref{teomain2}}
Define the constants $\alpha_{p}^{\L} = \frac{1}{2}\ELR(\llvert  \L\rrvert  ^p +
\llvert  \tL\rrvert  ^p)$
and $\alpha_{p}^{\R} = \frac{1}{2} \ELR(\llvert  \R \rrvert  ^p + \llvert  \tR\rrvert  ^p)$, so
$\alpha_{p} =
1-\alpha_{p}^{\L}-\alpha_{p}^{\R}$.
We first treat the case $p=1$.
Thus, we work with the processes $U^i$ solution of (\ref{eqUtk}) using
the functions $\G_t^i$ of Lemma~\ref{lemcoupling} with $p=1$.
Let us prove (i) first. We estimate the quantity
$f_t = \EE\llvert  X_t^1 - U_t^1\rrvert  $
which provides an upper bound for $\W_1(\law^1(\mathbf{X}_t),P_t)$.
Using (\ref{eqXtk}) and (\ref{eqUtk}), for all $0\leq s \leq t$ we have
%
\begin{eqnarray}
\label{eqXtUt} && \bigl\llvert X_t^1 -
U_t^1\bigr\rrvert - \bigl\llvert X_s^1
- U_s^1\bigr\rrvert\nonumber
\\
&&\qquad = \int_{(s,t]} \int_{\RR^2}\int
_{[0,N)} \bigl( \bigl\llvert \xi\bigl(X_{r^-}^1
- U_{r^-}^1\bigr)
+ \zeta\bigl(X_{r^-}^{\i(\tau)}- \G _r^1(\mathbf{X}_{r^-},\tau)\bigr)\bigr\rrvert
\nonumber\\[-8pt]\\[-8pt]\nonumber
&&\hspace*{66pt}\hspace*{151pt}{}  - \bigl\llvert X_{r^-}^1 - U_{r^-}^1
\bigr\rrvert \bigr)
\\
&&\hspace*{97pt}{}\times \N^1(dr,d\xi,d\zeta,d\tau).\nonumber
\end{eqnarray}
Recall that the intensity of $\N^1$ is $(N-1)^{-1} \,dt\barlr(d\xi,d\zeta)
\,d\tau\ind_{A^1}(\tau)$,
where $\barlr= (\law(\L,\R) + \law(\tL,\tR))/2$.
By the compensation formula, $t \mapsto f_t$ is absolutely continuous
and we obtain
%
\begin{eqnarray} \label{eqEXtUt}
f_t - f_s & \leq& \EE\int_s^t
\int_{\RR^2} \int_{A^1} \bigl( \bigl(
\llvert \xi \rrvert -1\bigr)\bigl\llvert X_r^1 -
U_r^1\bigr\rrvert + \llvert \zeta\rrvert \bigl\llvert
X_r^{\i(\tau)} - \G _r^1(\mathbf
{X}_r,\tau)\bigr\rrvert \bigr)\nonumber
\\
&&\hspace*{53pt}{}\times  \frac{d\tau}{N-1} \barlr(d\xi,d
\zeta) \,dr
\\
& =& \EE\int_s^t \bigl( \bigl(
\alpha_{1}^{\L} - 1\bigr)\bigl\llvert X_r^1
- U_r^1\bigr\rrvert + \alpha_{1}^{\R}
\W _1\bigl(\bar{\mathbf{X}}^{1}_r,P_r
\bigr) \bigr) \,dr,\nonumber
\end{eqnarray}
where in the last step we have used the fact that when $\tau$ is
uniform in $A^1$,
$(\G_s^1(\mathbf{x},\tau), x^{\i(\tau)})$ is an optimal coupling
between $P_s$ and $\bar{\mathbf{x}}^{1}$.
We deduce that for almost all $t\geq0$
%
\begin{equation}
\label{eqdft1} \partial_t f_t \leq-\bigl(1-
\alpha_{1}^{\L}\bigr) f_t +
\alpha_{1}^{\R} \EE \W_1\bigl(\bar{\mathbf
{X}}^{1}_t,P_t\bigr).
\end{equation}
Recall that $\bar{\mathbf{U}}^{i}_t = \frac{1}{N-1} \sum_{j\neq i}
\delta
_{U_t^j}$
for $i=1,\ldots,N$. The triangle inequality for $\W_1$ gives us
%
\begin{eqnarray} \label{eqW1XPt}
\EE\W_1\bigl(\bar{\mathbf{X}}^{1}_t,P_t
\bigr) &\leq&\EE\W_1\bigl(\bar{\mathbf{X}}^{1}_t,
\bar{\mathbf{U}}^{1}_t\bigr) + \EE \W _1\bigl(
\bar{\mathbf{U}}^{1}_t,P_t\bigr)
\nonumber\\[-8pt]\\[-8pt]\nonumber
&\leq&\EE\bigl\llvert X_t^1 - U_t^1
\bigr\rrvert + \EE\W_1\bigl(\bar{\mathbf{U}}^{1}_t,P_t
\bigr),
\end{eqnarray}
where the last inequality comes from the fact that $(X_t^{\i(\tau
)},U_t^{\i(\tau)})$
is a coupling between $\bar{\mathbf{X}}^{1}_t$ and $\bar{\mathbf
{U}}^{1}_t$ when $\tau$ is
uniformly chosen in $A^1$, and from the exchangeability
of $(X^i,U^i)_{i=1,\ldots,N}$. Putting this together with (\ref{eqdft1}),
we obtain
%
\begin{equation}
\label{eqEXU} \partial_t f_t \leq-
\alpha_{1} f_t + \alpha_{1}^{\R} \EE
\W_1\bigl(\bar{\mathbf{U}}^{1}_t,P_t
\bigr).
\end{equation}
Next, we need an estimate for $\EE\W_1(\bar{\mathbf{U}}^{1}_t,P_t)$.
Since the system
$(U^2,\ldots,U^N)$ is exchangeable, using a recent result
[Theorem 1.2 of \citet{hauray-mischler2014}], we obtain the
following: for each $q>0$ and each $\gamma< (2+1/q)^{-1}$, there
exists a constant $C_{q,\gamma}$ such that
%
\begin{equation}
\label{eqH&M} \qquad\EE\W_1\bigl(\bar{\mathbf{U}}^{1}_t,P_t
\bigr) \leq C_{q,\gamma} \mom_q(P_t)^{1/q}
\biggl( \W_1\bigl( \law\bigl(U_t^2,U_t^3
\bigr), P_t^{\otimes2}\bigr) + \frac{1}{N-1}
\biggr)^\gamma.
\end{equation}
Now, Lemma \ref{lemdecoupling} in the case $p=1$ and $k=2$ implies
$\W_1( \law(U_t^1,U_t^2), P_t^{\otimes2}) \leq C/N$,
where $C$ is some constant, which can change from line to line in what follows.
From this, Lemma \ref{lemmomentbounds}, and (\ref{eqEXU})--(\ref
{eqH&M}) we have
$\partial_t f_t \leq- \alpha_{1} f_t + C N^{-\gamma} e^{-(1/q)\barcq
{1,q} t}$,
and then Gronwall's lemma yields
\[
f_t \leq\frac{C}{N^\gamma} \int_0^t
e^{-\alpha_{1}(t-s)} e^{-(1/q)\barcq
{1,q} s}\,ds,
\]
since $f_0 = 0$.
Bounding $e^{-\alpha_{1} (t-s)} \leq e^{-(1/q)\barcq{1,q} (t-s)}$
gives (i) in the case $p=1$ and $k=1$.
From this and Lemma \ref{lemdecoupling}, case $k\geq2$ follows:
%
\begin{eqnarray}\label{eqcasekgeq2}
\W_1\bigl(\law^k(\mathbf{X}_t),P_t^{\otimes k}
\bigr)
&\leq&\W_1\bigl(\law^k(
\mathbf{X}_t),\law^k(\mathbf{U}_t) \bigr) +
\W_1\bigl(\law^k(\mathbf{U}_t),P_t^{\otimes k}
\bigr)
\nonumber\\[-8pt]\\[-8pt]\nonumber
&\leq&\EE\bigl\llvert X_t^1 - U_t^1
\bigr\rrvert + \frac{C k \min(1,t) e^{-\alpha_{1} t}}{N}.
\end{eqnarray}

We now prove (ii): as in (\ref{eqW1XPt}) we have
\begin{eqnarray*}
\EE\W_1(\bar{\mathbf{X}}_t,P_t) &\leq&\EE
\bigl\llvert X_t^1 - U_t^1\bigr
\rrvert + \EE\W_1(\bar{\mathbf{U}}_t,P_t)
\\
&\leq&\frac{Ct e^{-(1/q)\barcq{1,q} t }}{N^\gamma} + \frac{C
e^{-(1/q) \barcq{1,q} t}}{N^\gamma},
\end{eqnarray*}
where the last inequality comes from (i) in the case
$k=1$, and
from (\ref{eqH&M}) (with $\bar{\mathbf{U}}_t$ and $N$ in place of
$\bar{\mathbf{U}}^{1}_t$ and $N-1$)
together with Lemma \ref{lemdecoupling} in the case $k=2$.
From the previous inequality, (ii) follows;
moreover, the same estimate is also valid for $\EE\W_1(\bar{\mathbf
{X}}^{1}_t,P_t)$.

Now we treat the case $p=2$.
The proof is similar to the previous case,
with adaptations where required.
We work with the processes $U^i$ solution of (\ref{eqUtk}) using
the functions $\G_t^i$ of Lemma \ref{lemcoupling} with $p=2$.
As before,\vspace*{1pt} to prove the case $k=1$ we want to estimate $f_t = \EE
(X^1_t-U^1_t)^2$.
We proceed as in (\ref{eqXtUt}):
from (\ref{eqXtk}) and (\ref{eqUtk}), we have for all $0\leq s \leq t$
%
\begin{eqnarray}\label{eqXU2}
& &\bigl(X_t^1 - U_t^1
\bigr)^2 - \bigl(X_s^1 -
U_s^1\bigr)^2 \nonumber
\\
&&\qquad = \int_{(s,t]} \int_{\RR^2}\int
_{[0,N)} \bigl( \bigl[\xi\bigl(X_{r^-}^1 -
U_{r^-}^1\bigr)+ \zeta\bigl(X_{r^-}^{\i(\tau)}
- \G _r^1(\mathbf{X}_{r^-},\tau)\bigr)
\bigr]^2\nonumber
\\
&&\hspace*{65pt}\hspace*{152pt}{} - \bigl[X_{r^-}^1 -
U_{r^-}^1\bigr]^2 \bigr)\nonumber
\\
&&\hspace*{96pt}{}\times  \N^1(dr,d\xi,d\zeta,d\tau)
\\
&&\qquad = \int_{(s,t]} \int_{\RR^2}\int
_{[0,N)} \bigl( \bigl[\xi^2 - 1\bigr]
\bigl(X_{r^-}^1 - U_{r^-}^1
\bigr)^2+ \zeta^2 \bigl(X_{r^-}^{\i(\tau)}
- \G_r^1(\mathbf{X}_{r^-},\tau)
\bigr)^2 \nonumber
\\
&&\hspace*{142pt}{} + 2 \xi\zeta\bigl(X_{r^-}^1 -
U_{r^-}^1\bigr) \bigl(X_{r^-}^{\i
(\tau)} -
\G_r^1(\mathbf{X}_{r^-},\tau)\bigr) \bigr)\nonumber
\\
&& \hspace*{96pt}{}\times \N^1(dr,d\xi,d\zeta,d\tau ).
\nonumber
\end{eqnarray}
Taking expectations, the last term in the integral vanishes thanks to
condition $\ELR(\L\R+\tL\tR) = 0$. As in (\ref{eqEXtUt})--(\ref
{eqdft1}), this yields
%
\begin{equation}
\label{eqdft2} \partial_t f_t \leq-\bigl(1-
\alpha_{2}^{\L}\bigr) f_t +
\alpha_{2}^{\R} \EE \W_2^2\bigl(
\bar{\mathbf {X}}^{1}_t,P_t\bigr).
\end{equation}
Defining $g_t = \EE\W_2^2(\bar{\mathbf{U}}^{1}_t,P_t)$ and using the
triangle inequality of $\W_2$ we have
%
\begin{eqnarray} \label{eqEW2XPt}
&& \EE\W_2^2\bigl(\bar{\mathbf{X}}^{1}_t,P_t
\bigr)\nonumber
\\
&&\qquad   \leq\EE\W_2^2\bigl(\bar{\mathbf{X}}^{1}_t,
\bar{\mathbf{U}}^{1}_t\bigr) + 2 \EE\W_2
\bigl(\bar{\mathbf{X}}^{1}_t,\bar{\mathbf{U}}^{1}_t
\bigr)\W_2\bigl(\bar {\mathbf {U}}^{1}_t,P_t
\bigr) + \EE\W_2^2\bigl(\bar{\mathbf{U}}^{1}_t,P_t
\bigr)
\\
&&\qquad \leq f_t + 2 f_t^{1/2} g_t^{1/2}
+ g_t,\nonumber
\end{eqnarray}
where in the last inequality the term $f_t$ is obtained with the same argument
as in~(\ref{eqW1XPt}), and the term $f_t^{1/2}g_t^{1/2}$ comes from
the Cauchy--Schwarz inequality.
From this and (\ref{eqdft2}), we obtain
\[
\partial_t f_t \leq-\alpha_{2}
f_t + 2 \alpha_{2}^{\R
}f_t^{1/2}g_t^{1/2}
+ \alpha_{2}^{\R}g_t.
\]
Using a version of Gronwall's lemma [see, e.g.,  Lemma 4.1.8 of
\citet{ambrosio-gigli-savare2008}]
together with Jensen's inequality, we obtain
%
\begin{equation}
\label{eqEXU2} f_t \leq\alpha_{2}^{\R}
e^{-\alpha_{2} t} \bigl(2+8\alpha_{2}^{\R} t\bigr) \int
_0^t e^{\alpha_{2}
s} g_s \,ds.
\end{equation}
Now, we need an estimate for $g_t = \EE\W_2^2(\bar{\mathbf{U}}^{1}_t,P_t)$.
Unfortunately, we do not have at our disposal a result similar to
(\ref{eqH&M}),
which is valid only for $\W_1$.
To bypass this, we will make use of the following lemma (proved in
Section~\ref{secproofs}); it has the spirit of~(\ref{eqH&M})
in the sense that it will allow us to work with $\W_2^2(\law
^n(\mathbf
{U}_t), P_t^{\otimes n})$
instead of $\EE\W_2^2(\bar{\mathbf{U}}^{1}_t,P_t)$, but at the
price of
the extra term $\eps_{n,2}(P_t)$.

%
\begin{lem}
\label{lemW2Ymu}
Let $\mathbf{Y} = (Y^1,\ldots,Y^m)$ be an exchangeable random vector,
and let $\mu$ be a probability measure on $\RR$.
Then, for any $p\geq1$ and $n \leq m$, $n\in\NN$, we have
\begin{eqnarray*}
\frac{1}{2^{p-1}} \EE\W_p^p( \bar{\mathbf{Y}}, \mu)
&\leq& \frac{kn}{m} \bigl( \W_p^p\bigl(
\law^n(\mathbf{Y}),\mu^{\otimes
n}\bigr) + \eps_{n,p}(
\mu) \bigr)
\\
&&{} + \frac{\ell}{m} \bigl( \W_p^p\bigl(
\law^\ell(\mathbf{Y}),\mu ^{\otimes
\ell}\bigr) + \eps_{\ell,p}(
\mu) \bigr),
\end{eqnarray*}
where $k$ and $\ell$ are the unique nonnegative integers satisfying
$m = kn + \ell$, with $\ell\leq n-1$.
Here, $\eps_{n,p}(\mu):= \EE\W_p^p(\bar{\mathbf{Z}}, \mu)$,
where $\mathbf{Z}=(Z^1,\ldots,Z^n)$ are i.i.d. and $\mu$ distributed.
\end{lem}

Note that $\W_2^2(\law^\ell(\mathbf{U}_t),P_t^{\otimes\ell}) +
\eps
_{\ell,2}(P_t) \leq8 M_2(P_t)$.
Using this lemma with $p=2$, $m=N-1$, $\mathbf{Y} = (U^2_t,\ldots,U^{N-1}_t)$ and $\mu= P_t$, we obtain
that for every $n\leq N-1$
\begin{eqnarray*}
\EE\W_2^2\bigl(\bar{\mathbf{U}}^{1}_t,P_t
\bigr) &\leq& \W_2^2\bigl(\law^n(
\mathbf{U}_t),P_t^{\otimes n}\bigr) +
\eps_{n,p}(P_t) + \frac{n-1}{N-1} 8 M_2(P_t)
\\
&\leq& C \biggl( \frac{ne^{-\alpha_{2} t}}{N} + \eps_{n,2}(P_t)
\biggr),
\end{eqnarray*}
where in the last inequality we have used Lemmas \ref{lemmomentbounds}~and~\ref{lemdecoupling} with $p=2$ and $k=n$;
again $C$ is some constant that can change from line to line.
Putting this into (\ref{eqEXU2}) gives
%
\[
f_t \leq C (1+t) \biggl(\frac{nte^{-\alpha_{2} t}}{N} + \int
_0^t e^{-\alpha_{2}(t-s)} \eps_{n,2}(P_s)
\,ds \biggr).
\]
%
Given $q \in(2,q^*)$, $q\neq4$, from Theorem 1 of \citet
{fournier-guillin2013}
we know that $\eps_{n,2}(P_t) \leq C M_q^{2/q}(P_t) n^{-\eta}$, where
$\eta= \min(1/2,\frac{q-2}{q})$.
Choosing $n = \lfloor N^{1/(1+\eta)} \rfloor$ and using Lemma \ref
{lemmomentbounds}
with $p=2$ yields
%
\[
f_t \leq C (1+t) \biggl(\frac{te^{-\alpha_{2} t}}{N^\gamma} + \frac{1}{N^\gamma}\int
_0^t e^{-\alpha_{2}(t-s)} e^{-(2/q)
\barcq
{2,q} s} \,ds
\biggr),
\]
%
where\vspace*{1pt} $\gamma= \eta/(1+\eta) = \min(1/3, \frac{q-2}{2q-2})$.
Bounding $e^{-\alpha_{2}(t-s)} \leq e^{-(2/q) \barcq{2,q} (t-s)}$
gives (i) in the case $p=2$ and $k=1$.
The case $k\geq2$ follows as in (\ref{eqcasekgeq2}).

Finally, (ii) in the case $p=2$ follows from (\ref
{eqEW2XPt})
with a similar argument as in the case $p=1$. This completes the proof.
\end{pf*}

%
\begin{cor} \label{cortrajectorial}
Under the same hypotheses and notation of Theorem \ref{teomain2},
we have for all $T\geq0$,
\[
\EE\sup_{t\in[0,T]} \bigl\llvert X_t^1
- U_t^1\bigr\rrvert ^p \leq
\frac{C}{N^\gamma} \int_0^T
(1+t)^p e^{-(p/q)\barcq
{p,q}t } \,dt.
\]
\end{cor}

\begin{pf}
From (\ref{eqXtUt}), discarding
the negative term in the integral, we have
\begin{eqnarray*}
&& \sup_{t\in[0,T]} \bigl\llvert X_t^1 -
U_t^1\bigr\rrvert
\\
&&\qquad \leq\int_{(0,T]}
\int_{\RR^2}\int_{[0,N)} \bigl( \llvert \xi
\rrvert \bigl\llvert X_{t^-}^1 - U_{t^-}^1
\bigr\rrvert + \llvert \zeta\rrvert \bigl\llvert X_{t^-}^{\i(\tau)}
- \G _t^1(\mathbf {X}_{t^-},\tau)\bigr\rrvert
\bigr)
\\
&&\hspace*{99pt}{}\times  \N^1(dt,d\xi,d\zeta,d\tau).
\end{eqnarray*}
With the same argument that produced the term $\W_1(\bar{\mathbf
{X}}^{1}_r,P_r)$ in (\ref{eqEXtUt}),
the conclusion follows taking expectations and using
the previous estimates for $\EE\llvert  X_t^1 - U_t^1\rrvert  $ and $\EE\W_1(\bar
{\mathbf{X}}^{1}_t,P_t)$.
This proves the case $p=1$,
and the case $p=2$ follows from (\ref{eqXU2}) with a similar argument.
\end{pf}

%
\begin{rmk} \label{rmknon-cutoff}
To illustrate how our methods can indeed be used in noncutoff
contexts, consider Kac's model: $\L= \cos\theta= \tL$ and $\R=
-\sin\theta= -\tR$,
where $\theta$ is chosen according to an even cross-section function
$\beta\dvtx [-\pi,\pi] \to\RR_+$
that possibly is singular at $0$, but
satisfies the classical condition $\int_0^\pi\theta^2 \beta(\theta
)\,d\theta< \infty$,
see \citet{desvillettes-graham-meleard1999} for details.
Define $\beta_\eps(\theta) = \ind_{\llvert  \theta\rrvert  >\eps} \beta(\theta
)$ for a
given cutoff level $\eps>0$,
and associate with it the collection $(P_t^\eps)_{t\geq0}$ solving
$\partial_t P_t^\eps= \kappa_\eps(-P_t^\eps+ \Q^+_\eps
(P_t^\varepsilon))$,
where $\kappa_\eps= \int_{-\pi}^\pi\beta_\eps(\theta) \,d\theta$ and
$\Q_\eps^+$ is defined as
\begin{eqnarray*}
&& \int\phi(u) \Q_\eps^+(\mu) (du)
\\
&&\qquad = \int_\RR\int
_\RR\int_{-\pi
}^\pi\phi (u
\cos\theta- v \sin\theta) \frac{\beta_\eps(\theta)\,d\theta
}{\kappa
_\eps}\mu(dv) \mu(du).
\end{eqnarray*}
The particle system $\mathbf{X}^\eps$ and 
nonlinear processes $\mathbf{U}^\eps$ 
are constructed
in a way similar as in (\ref{eqXtk}) and (\ref{eqUtk})
but now using a Poisson measure $\N^{\eps,i}(dt,d\theta,d\tau)$
with intensity
$(N-1)^{-1}\,dt\,\beta_\eps(\theta)\,d\theta \,d\tau\ind_{A^i}(\tau)$
and functions $\G_t^{\eps,i}$ that couple optimally with $P_t^\eps$
instead of $P_t$.
Note that:
\begin{itemize}
\item The even moments of $P_t^\eps$ are controlled uniformly
in time and independently of $\eps$ [see, e.g., Lemma A.5 in
\citet{fournier-godinho2014}
in the case of the noncutoff nonlinear process; also, an induction similar
to the one used in the proof of Lemma \ref{lemmomentbounds} yields the
desired uniform bounds for $M_q(P_t^\eps)$ when $q$ is even].

\item The decoupling property of Lemma \ref{lemdecoupling} is also
valid for the system $\mathbf{U}_t^\eps$,
with constants independent of $\eps$:
in (\ref{eqEUV2}), all the terms involve either $1-L$ or $R^2$, which
correspond to $1-\cos\theta$ and $\sin^2 \theta$, respectively,
both of
order $\theta^2$.

\item In (\ref{eqEXU2}), the constant $\alpha_{2}^{\R}$ corresponds
to $\int_{-\pi}^\pi\sin^2 \theta\beta_\eps(\theta)$.
\end{itemize}

Thus, the argument can be replicated and the final constant will depend
on $\int_0^\pi\theta^2 \beta_\eps(\theta) \,d\theta$, which
remains bounded
as we let the cutoff $\eps\to0$.
Assuming, for instance, that $M_6(P_0) < \infty$,
this yields a constant $C$ independent of $\eps>0$ such that
\[
\EE\W_2^2\bigl(\bar{\mathbf{X}}_t^\eps,
P_t^\eps\bigr) \leq\frac{C(1+t)^2
}{ N^{1/3}}.
\]
%

However, we have not been able to obtain a trajectorial
result in the noncutoff case: discarding the negative term in the
integral of (\ref{eqXU2})
produces the term $\int_0^\pi\cos^2\theta\beta_\eps(\theta)
\,d\theta$
which no
longer stays bounded when $\eps\to0$.
\end{rmk}

\section{Proof of intermediate lemmas}
\label{secproofs}

\mbox{}

\begin{pf*}{Proof of Lemma \ref{lemcoupling}}
For fixed $n\in\NN$, given $\mathbf{y} = (y^1,\ldots,y^n) \in\RR^n$
recall that we write $\bar{\mathbf{y}} = \frac{1}{n} \sum_j \delta_{y^j}$.
The mapping $(t,\mathbf{y}) \mapsto(P_t,\bar{\mathbf{y}})$
from $\RR_+ \times\RR^n$ to $\P(\RR)\times\P(\RR)$ is
continuous when
$\P(\RR)$ is endowed with the weak topology (weak continuity of
$t\mapsto P_t$ is clear from the pathwise properties of the nonlinear
process). Thus, thanks to a measurable selection
result [see, e.g., Corollary~5.22 of \citet{villani2009}],
there exists a measurable mapping $(t,\mathbf{y}) \mapsto\pi_{t,\bar
{\mathbf{y}}}$ such that
$\pi_{t,\bar{\mathbf{y}}} \in\P(\RR\times\RR)$ is an optimal
transference plan
between $P_t$ and $\bar{\mathbf{y}}$.
We now define
\[
G(t,\mathbf{y},B) = \frac{\pi_{t,\bar{\mathbf{y}}}(B\times\{y^1\})}{\pi_{t,\bar
{\mathbf {y}}}(\RR\times\{y^1\})} = \pi_{t,\bar{\mathbf{y}}}\bigl(B\times\bigl
\{y^1\bigr\} \mid\RR\times\bigl\{y^1\bigr\} \bigr),
\]
for $t\geq0$, $\mathbf{y}\in\RR^n$ and any Borel set $B \subseteq
\RR$.
We claim that $G$ is a probability kernel from $\RR_+ \times\RR^n$
into $\RR$.
Indeed, it suffices to show that for every such $B$
the mapping $(t,\mathbf{y}) \mapsto\pi_{t,\bar{\mathbf
{y}}}(B\times
\{
y^1\})$ is measurable,
which in turn follows from the measurability of $(t,\mathbf{y})
\mapsto(\pi_{t,\bar{\mathbf{y}}},\mathbf{y})$ and the identity
\[
\pi_{t,\bar{\mathbf{y}}}\bigl(B\times\bigl\{y^1\bigr\}\bigr) = \lim
_{\eps\to0} \sum_{\ell\in\NN}
\pi_{t,\bar{\mathbf{y}}}\bigl(B \times D_\ell^{\eps}\bigr)
\ind_{D_\ell^{\eps}}\bigl(y^1\bigr),
\]
where $(D_\ell^{\eps})_{\ell\in\NN}$ is a measurable partition of
$\RR$
with $\operatorname{diam}(D_\ell^{\eps}) \leq\eps$.

Now, given $N\geq1$, with the kernel $G$ defined above for $n=N-1$ we
can associate a measurable mapping $g\dvtx \RR_+ \times\RR^{N-1}\times
[0,1]\to\RR$ or randomization of $G$
such that $g(t,\mathbf{y},\theta)$ has distribution $G(t,\mathbf
{y},\cdot)$
whenever $\theta$ is a uniform random variable in $[0,1]$
[see, e.g., Lemma 3.22 of \citet{kallenberg2002}]. For
$\mathbf{x}\in\RR^N$, we now put
%
\begin{equation}
\label{eqdefcoupfunct} \G_t^i(\mathbf{x},\tau) = \sum
_{j\neq i}^N \ind_{\{\i(\tau)=j\}} g\bigl(t,
\mathbf{x}^{(ij)}, \tau- \lfloor\tau\rfloor\bigr), \qquad\tau\in
A^i,
\end{equation}
where $\mathbf{x}^{(ij)} \in\RR^{N-1}$ denotes the vector $\mathbf{x}$
with its $i$ coordinate removed,
the $j$ coordinate in the first position, and the remaining
coordinates in positions $2,\ldots,N-1$
in increasing order.
We now show that when $\tau$ is uniform in $A^i$, $\G_t^i(\mathbf
{x},\tau)$ and $x^{\i(\tau)}$ have joint distribution
$\pi_{t,\bar{\mathbf{x}}^{i}}$.
Denoting $\IP^i$ the law of this
random variable $\tau$ and using the fact that $g(t,\mathbf
{x}^{(ij)},\theta)$ has law
$\pi_{t,\bar{\mathbf{x}}^{i}}(du \times\{ x^j \} \mid\RR\times\{
x^j\})$ when
$\theta$ is uniform in $[0,1]$,
we have for every fixed measurable set $B\subseteq\RR$ and every
$j\neq i$:
\begin{eqnarray*}
&& \IP^i\bigl(\G_t^i(\mathbf{x},\tau) \in B,
x^{\i(\tau)} = x^j\bigr)
\\
&&\qquad = \sum_{\ell\dvtx  x^\ell= x^j, \ell\neq i}
\int_{\ell-1}^\ell \ind_B\bigl(g
\bigl(t,\mathbf{x}^{(i\ell)},\tau- \lfloor\tau\rfloor\bigr)\bigr)
\frac
{d\tau}{N-1}
\\
&&\qquad = \frac{1}{N-1} \sum_{\ell\dvtx  x^\ell= x^j, \ell\neq i} \frac{\pi_{t,\bar{\mathbf{x}}^{i}}(B\times\{x^\ell\})}{\pi
_{t,\bar{\mathbf{x}}^{i}}(\RR\times\{x^\ell\})}
\\
&&\qquad = \frac{\llvert  \{\ell\dvtx  x^\ell= x^j, \ell\neq i \}\rrvert  }{(N-1)\pi_{t,\bar
{\mathbf{x}}^{i}}(\RR\times\{x^j\})} \pi_{t,\bar{\mathbf
{x}}^{i}}\bigl(B\times \bigl\{x^j
\bigr\}\bigr),
\end{eqnarray*}
where the quotient in the last line equals 1. This shows that $(\G
_t^i(\mathbf{x},\tau),x^{\i(\tau)})$ has distribution $\pi_{t,\bar
{\mathbf{x}}^{i}}$
and completes the proof of the existence of $\G^i$.

It remains to show that $\EE\int_{j-1}^j \phi(\G_t^i(\mathbf
{Y},\tau
)) \,d\tau= \langle P_t, \phi\rangle$
when $\mathbf{Y}$ is exchangeable, $j\neq i$ and $\phi$ is bounded and
measurable.
We get from (\ref{eqdefcoupfunct}) that
\begin{eqnarray*}
\int_{j-1}^j \phi\bigl(\G_t^i(
\mathbf{Y},\tau)\bigr) \,d\tau &=& \int_0^1 \phi
\bigl( g\bigl(t,\mathbf{Y}^{(ij)}, \tau\bigr)\bigr) \,d\tau
\\
&=& \int_\RR\phi(u) \pi_{t,\bar{\mathbf{Y}}^{i}}\bigl(du \times\bigl
\{Y^j\bigr\} \mid \RR\times\bigl\{Y^j\bigr\}\bigr),
\end{eqnarray*}
where we have again used that $g(t,\mathbf{Y}^{(ij)},\theta)$ has distribution
$\pi_{t,\bar{\mathbf{Y}}^{i}}(du \times\{Y^j\} \mid\RR\times\{
Y^j\}
)$ when $\theta$ is uniform in $[0,1]$.
From the exchangeability of $\mathbf{Y}$, it is clear that the last
expression has
the same distribution, for all $j\neq i$. Thus, its expected value
must be the
same for all $j\neq i$, and since
\[
\langle P_t, \phi\rangle = \int_{A^i} \phi
\bigl(\G_t^i(\mathbf{Y},\tau)\bigr) \frac{d\tau}{N-1} =
\sum_{j\neq i} \int_{j-1}^j
\phi\bigl(\G_t^i(\mathbf{Y},\tau)\bigr)
\frac
{d\tau}{N-1},
\]
the conclusion follows.
\end{pf*}

%
\begin{rmk}\label{rmkpostthmcoupling}
Since we are working on $\RR$, the increasing coupling between $P_t$
and $\bar{\mathbf{x}}^{i}$
is in fact an optimal coupling [see, e.g., Theorem 6.0.2 in
\citet{ambrosio-gigli-savare2008}],
which allows for a simpler proof of Lemma \ref{lemcoupling}. However,
we opted to give a proof
that remains valid on $\RR^d$ with the hope that this coupling can be
used in a more general setting.
\end{rmk}

\begin{pf*}{Proof of Lemma \ref{lemXU}}
Existence and uniqueness for (\ref{eqUtk}) are obtained with a
construction similar to (\ref{eqXtn}).
To show that $U^i$ is a nonlinear process,
define $\tilde{\N}^i(dt,d\xi,d\zeta,dv)$ to be the point measure
on $\RR_+
\times\RR^2 \times\RR$
with atoms $(t,\xi,\zeta,\G^i_t(\mathbf{X}_{t^-},\tau))$ for
every atom
$(t,\xi,\zeta,\tau)$ of $\N^i$;
since the dependence on $\mathbf{X}$ is predictable, one can use the
compensation formula to compute the Laplace functional of $\tilde{\N}^i$
and conclude that $\tilde{\N}^i$ is a Poisson point measure with intensity
$dt \barlr(d\xi,d\zeta) P_t(dv)$.
Then (\ref{eqnlpr}) is satisfied for $V=U^i$ with
$\M= \tilde{\N}^i$, implying that $\law(U^i) = P$. The collection
$(X^1,U^1),\ldots,(X^N,U^N)$ is obviously exchangeable.
\end{pf*}

\begin{pf*}{Proof of Lemma \ref{lemmomentbounds}}
Call $h_t^q = \int\llvert  u\rrvert  ^q P_t(du)$.
We first prove the statement for the case $p=2$.
Using (\ref{eqPt})--(\ref{eqQ+}) with $\phi= \llvert  \cdot\rrvert  ^2$
yields $\partial_t h_t^2 = - \alpha_{2} h_t^2 + \ELR(\L\R+ \tL\tR)
(h_t^1)^2$,
and since $\ELR(\L\R+ \tL\tR) = 0$ this implies $h_t^2 = h_0^2
e^{-\alpha_{2}t}$.
Assume now that $q \in(2,q^*)$ is an integer.
Using (\ref{eqPt})--(\ref{eqQ+}) with $\phi= \llvert  \cdot\rrvert  ^q$, we have
%
\begin{eqnarray}\label{eqhtq}
\partial_t h_t^q &=& - h_t^q
+ \frac{1}{2} \iint\ELR\bigl( \llvert \L u + \R v\rrvert ^q +
\llvert \tL v + \tR u\rrvert ^q \bigr) P_t(du)
P_t(dv)
\nonumber\\[-8pt]\\[-8pt]\nonumber
&\leq& - \alpha_{q} h_t^q + \frac{1}{2}
\sum_{i=1}^{q-1} \pmatrix{{q}
\vspace*{3pt}\cr
{i}}
h_t^i h_t^{q-i} \ELR\bigl( \llvert
\L \rrvert ^i\llvert \R \rrvert ^{q-i} + \llvert \tL\rrvert
^i\llvert \tR\rrvert ^{q-i} \bigr).
\nonumber
\end{eqnarray}
Using loose bounds for $\binomt{q}{i}$, we obtain
\[
h_t^q \leq h_0^q
e^{-\alpha_{q}t} + C \sum_{i=1}^{q-1} \int
_0^t e^{ -\alpha_{q}(t-s)} h_s^i
h_s^{q-i} \,ds,
\]
where $C$ is a constant that does not depend on $t$, and may change
from line to line.
We now apply induction: the case $q=2$ was already proven, and for
$q\in
(2,q^*)$ integer,
assuming the desired property for all integer in $\{2,\ldots,q-1\}$ and
using the bound $h_t^1 \leq(h_t^2)^{1/2} \leq C e^{-\alpha_{2}t/2}$,
we obtain
\begin{eqnarray*}
h_t^q &\leq& h_0^q
e^{-\alpha_{q}t} + C \int_0^t
e^{ -\alpha_{q}(t-s)} e^{-\alpha_{2}s/2} e^{-\barcq
{2,q-1} s} \,ds
\\
&&{} + C \sum_{i=2}^{q-2} \int
_0^t e^{ -\alpha_{q}(t-s)} e^{-\barcq
{2,i} s}
e^{-\barcq{2,q-i} s} \,ds.
\end{eqnarray*}
Note that $\alpha_{q}>0$, since $2<q<q^*$, and recall that
$\barcq{2,q}:= \inf_{2\leq r \leq q} \alpha_{r} = \min(\alpha
_{2},\alpha_{q})$. Thus,
if $\alpha_{2} = 0$ then $\barcq{2,i} = \barcq{2,q-i} = \barcq{2,q} =
0$ and
the last inequality yields
$h_t^q \leq h_0^q + C \int_0^t e^{-\alpha_{q}(t-s)}\,ds \leq C$, as desired.
On the other hand, if $\alpha_{2} > 0$, we bound $\alpha_{2}$,
$\alpha_{q}$,
$\barcq
{2,i}$ and $\barcq{2,q-i}$
from\vspace*{1pt} below by $\barcq{2,q}>0$ and obtain
$h_t^q \leq h_0^q e^{-\barcq{2,q}t} + C \int_0^t e^{-\barcq
{2,q}(t+s/2)}\,ds \leq C e^{-\barcq{2,q}t}$,
which completes the induction and the proof in the case $p=2$ and
integer $q\in\{2\} \cup(2,q^*)$.

Assume now that $2<q=m+\eps<q^*$ with $m\in\{2,\ldots\}$ and $\eps
\in(0,1)$.
Bounding $\llvert  x+y\rrvert  ^{q} \leq(\llvert  x\rrvert  +\llvert  y\rrvert  )^m(\llvert  x\rrvert  ^\eps+\llvert  y\rrvert  ^\eps)$ in (\ref{eqhtq})
and using the binomial theorem as before, we obtain
\begin{eqnarray*}
\partial_t h_t^q &=& - h_t^q
+ \frac{1}{2} \sum_{i=0}^m
\pmatrix{{m}
\vspace*{3pt}\cr
{i}} \iint \ELR \bigl( \bigl(\llvert \L u\rrvert ^\eps+
\llvert \R v\rrvert ^\eps\bigr) \llvert \L u\rrvert ^i
\llvert \R v\rrvert ^{m-i}
\\
&&\hspace*{112pt}{} + \bigl(\llvert \tL v\rrvert ^\eps+
\llvert \tR u\rrvert ^\eps\bigr) \llvert \tL v\rrvert ^i
\llvert \tR u\rrvert ^{m-i} \bigr) P_t(du)
P_t(dv)
\\
&=& -\alpha_{q} h_t^q + \sum
_{i=0}^{m-1} \pmatrix{{m}
\vspace*{3pt}\cr
{i}} \ELR \bigl(\llvert \L
\rrvert ^{i+\eps}\llvert \R \rrvert ^{m-i} + \llvert \L\rrvert
^{m-i}\llvert \R \rrvert ^{i+\eps}
\\
&&\hspace*{100pt}{} + \llvert \tL\rrvert
^{i+\eps}\llvert \tR\rrvert ^{m-i} + \llvert \tL\rrvert
^{m-i}\llvert \tR\rrvert ^{i+\eps} \bigr) h_t^{i+\eps}h_t^{m-i},
\end{eqnarray*}
which yields
\[
h_t^q \leq h_0^q
e^{-\alpha_{q}t} + C \sum_{i=0}^{m-1} \int
_0^t e^{ -\alpha_{q}(t-s)} h_s^{i+\eps
}h_s^{m-i}
\,ds.
\]
Note that $h_t^r \leq(h_t^2)^{r/2} \leq C e^{-r \alpha_{2} t /2}$ for $r
\in(0,2)$.
This and the fact that the property is true for the integers,
allow us to use induction on $m$ in a way similar as before, and
complete the proof in the case $p=2$.

A similar argument, with the induction starting at $q=1$, proves the
case $p=1$.
\end{pf*}

\begin{pf*}{Proof of Lemma \ref{lemdecoupling}}
Let us first prove the case $p=1$.
Given $k\in\{2,\ldots,N\}$ fixed, we want to construct $k$ independent
nonlinear processes
$V^1,\ldots,V^k$ such that $\EE\llvert  U_t^i - V_t^i\rrvert  $ is small. To achieve
this, we will decouple
$U^1,\ldots,U^k$ by replacing the shared atoms of $\N^1,\ldots,\N^k$
with new, independent atoms. To this end, let $\M$ be an independent
copy of
$\N$ (also independent from $\mathbf{X}_0$), 
and define for each $i\in\{1,\ldots,k\}$
%
\begin{eqnarray}
\label{eqMkl}
&& \M^{i}(dt,d\xi,d\zeta,d\tau)\nonumber
\\
&&\qquad = \N\bigl(dt,\bigl(d\xi
\times d\zeta\times\RR^2\bigr),[i-1,i),d\tau\bigr)
\nonumber\\[-8pt]\\[-8pt]\nonumber
&&\quad\qquad{}+ \N\bigl(dt,\bigl(\RR^2 \times d\xi\times d\zeta\bigr),d
\tau,[i-1,i)\bigr) \ind _{[k,N)}(\tau)
\\
&&\quad\qquad{} + \M\bigl(dt,\bigl(\RR^2 \times d\xi\times d\zeta\bigr),d
\tau,[i-1,i)\bigr) \ind _{[0,k)}(\tau).\nonumber
\end{eqnarray}
Note that $\M^{i}$ is, like $\N^i$, a Poisson point
measure on $\RR_+ \times\RR^2 \times[0,N)$ with intensity
$(N-1)^{-1}\,dt \barlr(\xi,\zeta) \,d\tau\ind_{A^i}(\tau)$,
and that $\M^{1},\ldots,\M^{k}$ are independent.
Following (\ref{eqUtk}), we define $V^{i}$ as the solution of
%
\begin{equation}
\label{eqVtkl} dV_t^{i} = \int_{\RR^2}
\int_{[0,N)} \bigl[(\xi-1)V_{t^-}^{i} +
\zeta\G _t^i(\mathbf {X}_{t^-},\tau)\bigr]
\M^{i}(dt,d\xi,d\zeta,d\tau),
\end{equation}
with $V_0^{i} = U_0^i$.
If we define $\tilde{\M}^i$ to be the point process in $\RR_+ \times
\RR
^2 \times\RR$
with atoms $(t,\xi,\zeta,\G_t^i (\mathbf{X}_{t^-},\tau))$ for
every atom
$(t,\xi,\zeta,\tau)$ of $\M^i$,
it is clear that $V^i$ depends only on $\tilde{\M}^i$ and $X_0^i$.
Since: (i) the dependence on $\mathbf{X}$ is predictable,
(ii) the Poisson measures $\M^1,\ldots,\M^k$ are independent and
(iii) the $\tau$-law of $\G_t^i(\mathbf{x},\tau)$ is $P_t$ for every
$\mathbf{x}\in\RR^N$,
one can use the compensation formula to compute the joint
Laplace functional of $\tilde{\M}^1,\ldots,\tilde{\M}^k$
and conclude that they are independent Poisson point measures,
all with intensity $dt \barlr(d\xi,d\zeta) P_t(dv)$.
This shows that each $V^i$ is a nonlinear process and that
they are independent.

Consequently, we have
\[
\W_1\bigl( \law^k(\mathbf{U}_t),
P_t^{\otimes k}\bigr) \leq\EE \Biggl( \frac{1}{k} \sum
_{i=1}^k \bigl\llvert U_t^i
- V_t^{i}\bigr\rrvert \Biggr) = \EE\bigl\llvert
U_t^1 - V_t^{1}\bigr\rrvert,
\]
where in the last step we used the fact that all the $(U^i,V^{i})$'s
have the same law. To estimate the last term
$h_t = \EE\llvert  U_t^1 - V_t^{1}\rrvert  $, we proceed as in (\ref{eqXtUt}):
from (\ref{eqNk}), (\ref{eqUtk}), (\ref{eqMkl}) and (\ref{eqVtkl}),
we have
for all $0\leq s\leq t$:
%
\begin{equation}
h_t = h_s + \EE\int_{(s,t]} \int
_{\RR^2}\int_{[0,N)} \bigl(
J_r^1 + J_r^2 +
J_r^3 \bigr), \label{eqUtkVtkl}
\end{equation}
where $J_r^1$ is the term associated with the simultaneous jumps of
$U^1$ and $V^1$,
$J_r^2$ corresponds to the jumps of $U^1$ alone, and $J_r^3$ gives the
jumps of $V^1$ alone.
Specifically,
\begin{eqnarray*}
J_r^1 &=& \bigl(\llvert \xi\rrvert -1\bigr) \bigl\llvert
U_{r^-}^1 - V_{r^-}^{1}\bigr\rrvert
\bigl(  \N\bigl(dr,\bigl(d\xi\times d\zeta\times\RR^2
\bigr),[0,1),d\tau\bigr)
\\
&&\hspace*{93pt}{}  + \N\bigl(dr,\bigl(\RR^2 \times d\xi
\times d\zeta\bigr),d\tau,[0,1)\bigr) \ind_{[k,N)}(\tau) \bigr),
\\
J_r^2 &=& \bigl(\bigl\llvert \xi\bigl(U_{r^-}^1
- V_{r^-}^{1}\bigr) + \zeta\G_r^1(
\mathbf {X}_{r^-},\tau ) + (\xi-1)V_{r^-}^{1}\bigr
\rrvert - \bigl\llvert U_{r^-}^1 - V_{r^-}^{1}
\bigr\rrvert \bigr)
\\
&&{}\times \N\bigl(dr,\bigl(\RR^2 \times d\xi\times d\zeta
\bigr),d\tau,[0,1)\bigr) \ind_{[0,k)}(\tau),
\\
J_r^3 &=& \bigl( \bigl\llvert \xi\bigl(U_{r^-}^1
- V_{r^-}^{1}\bigr) - \zeta\G_r^1(
\mathbf {X}_{r^-},\tau ) - (\xi-1)U_{r^-}^1\bigr
\rrvert - \bigl\llvert U_{r^-}^1 - V_{r^-}^{1}
\bigr\rrvert \bigr)
\\
&&{}\times \M\bigl(dr,\bigl(\RR^2 \times d\xi\times d\zeta
\bigr),d\tau,[0,1)\bigr) \ind_{[0,k)}(\tau).
\end{eqnarray*}
Then
\[
\EE\int_{(s,t]} \int_{\RR^2} \int
_{[0,N)} J_r^1 = \biggl(
\frac{1}{2}\bigl(\ELR\llvert \L\rrvert -1\bigr) + \frac{1}{2}\bigl(
\ELR\llvert \tL\rrvert -1\bigr) \frac
{N-k}{N-1} \biggr) \int
_s^t h_r \,dr.
\]
Using the triangle inequality in the term $J_r^2$,
\begin{eqnarray*}
&& \EE\int_{(s,t]} \int_{\RR^2} \int
_{[0,N)} J_r^2
\\
&&\qquad \leq\EE\int
_s^t \int_1^k
\bigl( \bigl(\ELR\llvert \tL\rrvert -1\bigr) \bigl\llvert U_{r}^1
- V_{r}^{1}\bigr\rrvert + \ELR\llvert \tR\rrvert \bigl
\llvert \G_r^1(\mathbf{X}_r,\tau)\bigr
\rrvert + \ELR \llvert \tL -1\rrvert \bigl\llvert V_r^{1}
\bigr\rrvert \bigr)
\\
&&\hspace*{68pt}{}\times  \frac{dr \,d\tau}{2(N-1)}.
\end{eqnarray*}
From Lemma \ref{lemcoupling},
we know that $\EE\int_{i-1}^i \llvert  \G_r^1(\mathbf{X}_r,\tau)\rrvert   \,d\tau=
\mom
_1(P_r)$ for all $i=2,\ldots,k$.
Using that $V_r^{1}$ has law $P_r$, we obtain
\begin{eqnarray*}
&&\EE\int_{(s,t]} \int_{\RR^2} \int
_{[0,N)} J_r^2
\\
&&\qquad \leq\frac{k-1}{2(N-1)} \biggl( \bigl(\ELR\llvert \tL\rrvert -1\bigr) \int
_s^t h_r \,dr + \bigl(\ELR\llvert
\tL-1\rrvert + \ELR\llvert \tR\rrvert \bigr) \int_s^t
\mom_1(P_r) \,dr \biggr).
\end{eqnarray*}
With a similar argument, the last inequality is also valid with $J_r^3$
in the left-hand side.
Putting all this into (\ref{eqUtkVtkl}), we have
\begin{eqnarray*}
h_t &\leq& h_s - \biggl(1-\frac{1}{2}\ELR\bigl(
\llvert \L\rrvert +\llvert \tL\rrvert \bigr) + \frac
{1}{2}\bigl(1-\ELR
\llvert \tL\rrvert \bigr)\frac{k-1}{N-1} \biggr) \int_s^t
h_r \,dr
\\
&&{} + \frac{(\ELR\llvert  \tL-1\rrvert   + \ELR\llvert  \tR\rrvert  ) (k-1)}{N-1}\int_s^t
\mom_1(P_r)\,dr.
\end{eqnarray*}
Recall the constants $\alpha_{1}^{\L} = \frac{1}{2}\ELR(\llvert  \L\rrvert  +\llvert  \tL\rrvert  )$,
$\alpha_{1}^{\R} = \frac{1}{2}\ELR(\llvert  \R \rrvert  +\llvert  \tR\rrvert  )$ and $\alpha_{1}
= 1- \alpha_{1}^{\L}
- \alpha_{1}^{\R}$.
Also, put $\bi= \frac{1}{2}(1-\ELR\llvert  \tL\rrvert  )$, which can be assumed
nonnegative without loss of generality
[if not, exchange the roles of $(\L,\R)$ and $(\tL,\tR)$].
From the previous inequality and from Lemma \ref{lemmomentbounds} in
the case $q=1$,
it follows that for almost all $t\geq0$,
\[
\partial_t h_t \leq- \biggl( \alpha_{1} +
\alpha_{1}^{\R} + \bi \frac{k-1}{N-1} \biggr)
h_t + \frac{C (k-1) e^{-\alpha_{1} t} }{N-1},
\]
and now Gronwall's lemma gives
\[
h_t \leq\frac{C(k-1)e^{-\alpha_{1} t}}{(N-1)(\alpha_{1}^{\R} + \bi
((k-1)/(N-1)))} \bigl[ 1 - e^{-(\alpha_{1}^{\R} + \bi ((k-1)/(N-1)))t} \bigr].
\]
Using the inequality $1-e^{-x} \leq x$, the desired result follows for
the case $p=1$.

In\vspace*{1pt} the case $p=2$, we construct the system $V^1,\ldots,V^k$ exactly as before,
but using the functions $\G_t^i$ provided by Lemma \ref{lemcoupling}
with cost $\llvert  x-y\rrvert  ^2$.
To obtain the desired inequality for
$\W_2^2(\law^k(\mathbf{U}_t),P_t^{\otimes k})$, it suffices to work
with $h_t = \EE(U_t^1-V_t^1)^2$.
We also have (\ref{eqUtkVtkl}), where $J_r^1$, $J_r^2$ and $J_r^3$ now
are given by
\begin{eqnarray*}
J_r^1 &=& \bigl(\xi^2-1\bigr)
\bigl(U_{r^-}^1 - V_{r^-}^{1}
\bigr)^2 \bigl(  \N\bigl(dr,\bigl(d\xi\times d\zeta
\times\RR^2\bigr),[0,1),d\tau\bigr)
\\
&&\hspace*{99pt}{} + \N\bigl(dr,\bigl(
\RR^2 \times d\xi\times d\zeta\bigr),d\tau,[0,1)\bigr)
\ind_{[k,N)}(\tau) \bigr),
\\
J_r^2 &=& \bigl(\bigl(\xi\bigl(U_{r^-}^1
- V_{r^-}^{1}\bigr) + \zeta\G_r^1(
\mathbf {X}_{r^-},\tau ) + (\xi-1)V_{r^-}^{1}
\bigr)^2 - \bigl(U_{r^-}^1 -
V_{r^-}^{1}\bigr)^2 \bigr)
\\
&&{}\times \N\bigl(dr,\bigl(\RR^2 \times d\xi\times d\zeta
\bigr),d\tau,[0,1)\bigr) \ind_{[0,k)}(\tau),
\\
J_r^3 &=& \bigl( \bigl(\xi\bigl(U_{r^-}^1
- V_{r^-}^{1}\bigr) - \zeta\G_r^1(
\mathbf {X}_{r^-},\tau ) - (\xi-1)U_{r^-}^1
\bigr)^2 - \bigl(U_{r^-}^1 -
V_{r^-}^{1}\bigr)^2 \bigr)
\\
&&{}\times \M\bigl(dr,\bigl(\RR^2 \times d\xi\times d\zeta
\bigr),d\tau,[0,1)\bigr) \ind_{[0,k)}(\tau).
\end{eqnarray*}
Using that $\EE\int_{i-1}^i \G_t^1(\mathbf{X}_t,\tau)^2 \,d\tau=
\mom
_2(P_t)$ for all $i=2,\ldots,k$,
we obtain
\begin{eqnarray*}
&& \EE\int_{(s,t]} \int_{\RR^2} \int
_{[0,N)} J_r^1 = \biggl(
\frac{1}{2}\bigl(\ELR\L^2-1\bigr) + \frac{1}{2}\bigl(\ELR
\tL^2-1\bigr) \frac
{N-k}{N-1} \biggr) \int_s^t
h_r \,dr,
\\[-1pt]
&& \EE\int_{(s,t]} \int_{\RR^2} \int
_{[0,N)} J_r^2
\\[-1pt]
&&\qquad = \int
_s^t \biggl( \bigl(\ELR
\tL^2 -1\bigr) h_r + \ELR \tR^2
\mom_2(P_r) + \ELR(\tL-1)^2
\mom_2(P_r)
\\[-1pt]
&&\hspace*{52pt}{}+ 2 \ELR(\tL\tR) \EE\bigl(U_r^1-V_r^1
\bigr)\int_1^k \G_r^1(
\mathbf {X}_r,\tau ) \frac{d\tau}{k-1}
\\[-1pt]
&&\hspace*{52pt}{} + 2\ELR\bigl(\tL(\tL-1)\bigr)
\EE\bigl(U_r^1-V_r^1\bigr)
V_r^{1}
\\[-1pt]
&&\hspace*{32pt}\hspace*{52pt}{} + 2\ELR\bigl((\tL-1)\tR\bigr) \EE
V_r^1 \int_1^k
\G_r^1(\mathbf {X}_r,\tau)\frac{d\tau}{k-1}
\biggr) \frac{(k-1) \,dr}{2(N-1)},
\\[-1pt]
&& \EE\int_{(s,t]} \int_{\RR^2} \int
_{[0,N)} J_r^3
\\[-1pt]
&&\qquad = \int
_s^t \biggl( \bigl(\ELR
\tL^2 -1\bigr) h_r + \ELR \tR^2
\mom_2(P_r) + \ELR(\tL-1)^2
\mom_2(P_r)
\\[-1pt]
&&\hspace*{51pt}{} - 2 \ELR(\tL\tR) \EE\bigl(U_r^1-V_r^1
\bigr)\int_1^k \G_r^1(
\mathbf {X}_r,\tau ) \frac{d\tau}{k-1}
\\[-1pt]
&&\hspace*{51pt}{} - 2\ELR\bigl(\tL(\tL-1)\bigr)
\EE\bigl(U_r^1-V_r^1\bigr)
U_r^{1}
\\[-1pt]
&&\hspace*{32pt}\hspace*{51pt}{}+ 2\ELR\bigl((\tL-1)\tR\bigr) \EE
U_r^1 \int_1^k
\G_r^1(\mathbf {X}_r,\tau)\frac{d\tau}{k-1}
\biggr) \frac{(k-1) \,dr}{2(N-1)}.
\end{eqnarray*}
From this and (\ref{eqUtkVtkl}), we have for almost all $t\geq0$
%
\begin{eqnarray}
\label{eqEUV2} \partial_t h_t &=& -h_t
\biggl( \biggl(1-\frac{1}{2}\ELR\bigl(\L^2 + \tL^2
\bigr) \biggr) + \frac
{k-1}{2(N-1)}\ELR(\tL-1)^2 \biggr)\nonumber
\\[-1pt]
&&{} + \mom_2(P_t) \frac{k-1}{N-1} \bigl(
\ELR(\tL-1)^2 + \ELR\tR ^2\bigr)
\\[-1pt]
&&{} + \ELR\bigl((\tL-1)\tR\bigr) \EE\bigl(U_t^1
+ V_t^1\bigr) \int_1^k
\G _t^1(\mathbf{X}_t,\tau)
\frac{d\tau}{N-1}.\nonumber
\end{eqnarray}
We also have $\int_{i-1}^i \EE(U_t^1 + V_t^1) \G_t^1(\mathbf
{X}_t,\tau
) \,d\tau\leq2 \mom_2(P_t)$
for all $i=2,\ldots,k$, thanks to the Cauchy--Schwarz and Jensen
inequalities.
Recall the constants
$\alpha_{2}^{\L} = \frac{1}{2}\ELR(\L^2 + \tL^2)$, $\alpha
_{2}^{\R} = \frac
{1}{2}\ELR(\R
^2 + \tR^2)$,
$\alpha_{2} = 1-\alpha_{2}^{\L} - \alpha_{2}^{\R}$,
and put $\bii= \frac{1}{2}\ELR(\tL-1)^2$.
Using Lemma \ref{lemmomentbounds} in the case $p=q=2$,
we thus obtain
\[
\partial_t h_t \leq- \biggl(\alpha_{2} +
\alpha_{2}^{\R} + b\frac{k-1}{N-1} \biggr)h_t
+ \frac{C (k-1) \mom_2(P_0) e^{-\alpha_{2} t}}{N-1},
\]
and the conclusion follows from Gronwall's lemma as before.
\end{pf*}

\begin{pf*}{Proof of Lemma \ref{lemW2Ymu}}
For simplicity, we will prove only the case $\ell=0$,
that is, when $n$ divides $m$. Let us arrange a vector $\mathbf{y}\in
\RR
^m$ as a matrix
with $k$ rows and $n$ columns, that is, $\mathbf{y} = (y^{ij})$, with
$i=1,\ldots,k$,
$j=1,\ldots,n$, and write $\mathbf{y}_i = (y^{i1},\ldots,y^{in})$
and $\bar{\mathbf{y}}_i = \frac{1}{n}\sum_{j=1}^n \delta_{y^{ij}}$.
Let us couple $\mathbf{Y}$ with a random vector $\mathbf{Z} \in(\RR^n)^k$
in such a way that each $(\mathbf{Y}_i,\mathbf{Z}_i)$ is an optimal
coupling between
$\law^n(\mathbf{Y})$
and $\mu^{\otimes n}$
[with respect to the cost function $d_{n,p}^p(\cdot,\cdot)$ of (\ref
{eqdkp}), as usual].
Using the latter, we have
%
\begin{equation}
\label{eqEWpYZ} \EE\W_p^p(\bar{\mathbf{Y}},\bar{
\mathbf{Z}}) \leq\frac{1}{k} \sum_{i=1}^k
\EE\frac{1}{n} \sum_{j=1}^n \bigl
\llvert Y^{ij}-Z^{ij}\bigr\rrvert ^p =
\W_p^p\bigl(\nu^n,\mu^{\otimes n}\bigr).
\end{equation}
On the other hand, for each $i=1,\ldots,k$ there is a function $q^i\dvtx
\RR^m \times[0,1] \to\RR$
such that for all $\mathbf{z} \in\RR^m$, the pair $(z^{i\i^n(\theta
)},q^i(\mathbf{z},\theta))$ with $\i^n(\theta) = \i(n\theta) =
\lfloor
n \theta\rfloor+ 1$, is an optimal coupling
between $\bar{\mathbf{z}}_i$ and $\mu$ when $\theta$ is uniformly chosen
in $[0,1]$.
Now we randomize the choice of $i$ with a uniform variable $\vartheta
\in[0,1]$
independent of $\theta$, so $z^{\i^k(\vartheta) \i^n(\theta)}$ and
$q^{\i^k(\vartheta)}(\mathbf{z},\theta)$
are $(\theta,\vartheta)$-realizations of $\bar{\mathbf{z}}$ and
$\mu$,
respectively.
Putting $\mathbf{Z}$ in place of $\mathbf{z}$, this construction gives
\begin{eqnarray*}
\EE\W_p^p(\bar{\mathbf{Z}},\mu) &\leq&\EE
\iint_{[0,1]^2} \bigl\llvert Z^{\i^k(\vartheta) \i^n(\theta)} - q^{\i
^k(\vartheta)}(
\mathbf{Z},\theta)\bigr\rrvert ^p \,d\vartheta \,d\theta
\\
&=& \EE\frac{1}{k} \sum_{i=1}^k
\W_p^p(\bar{\mathbf{Z}}_i,\mu) =
\eps_{n,p}(\mu).
\end{eqnarray*}
[Recall that $\eps_{n,p}(\mu) = \EE\W_p^p(\frac{1}{n}\sum_i
\delta
_{\zeta^i},\mu)$,
with $\zeta^1,\ldots,\zeta^n$ independent and \mbox{$\mu$-}distributed].
With this and (\ref{eqEWpYZ}), we conclude in the case $\ell=0$
\begin{eqnarray*}
\EE\W_p^p(\bar{\mathbf{Y}},\mu) &\leq&\EE \bigl(
\W_p(\bar{\mathbf{Y}},\bar{\mathbf{Z}}) + \W _p(\bar{
\mathbf {Z}},\mu) \bigr)^p
\leq 2^{p-1} \bigl(\EE\W_p^p(\bar{\mathbf{Y}},
\bar{\mathbf{Z}}) + \EE\W _p^p(\bar{\mathbf{Z}},\mu)
\bigr).
\end{eqnarray*}
In the case $\ell>0$, the construction is similar, but now $(\mathbf
{Y},\mathbf{Z})$
must include an additional optimal coupling between
$\law^\ell(\mathbf{Y})$
and $\mu^{\otimes\ell}$,
which gives the extra term.
\end{pf*}


\section*{Acknowledgements}
We thank two anonymous referees for carefully reading a former version
of this work and for their questions and remarks that allowed us to
improve its presentation.


%

\printaddresses
\end{document}